\documentclass[12pt]{article}

\usepackage{amsfonts}
\usepackage{amsmath}
\usepackage{amsfonts}
\usepackage{amsmath}
\usepackage{epsfig}
\usepackage[compress]{cite}
\usepackage{graphics}
\usepackage[papersize={8.5in,10.5in},margin=1in]{geometry}

\begin{document}

\pagestyle{plain}

\begin{center}\Large
\bf{Implicit sampling for path integral control,\\Monte Carlo localization and SLAM}
\end{center}

\begin{center}
Matthias Morzfeld\\
Department of Mathematics\\
University of California, Berkeley,\\
and\\
Lawrence Berkeley National Laboratory
\end{center}

\begin{center}
\emph{Abstract}
\end{center}
Implicit sampling is a recently-developed variationally-enhanced sampling method,
that guides its samples to regions of high probability, so that each
sample carries information. 
Implicit sampling may thus improve the performance of algorithms that rely on Monte Carlo methods.
Here the applicability and usefulness of implicit sampling 
for improving the performance of Monte Carlo methods in 
estimation and control is explored, and
implicit sampling based algorithms for stochastic optimal control, 
stochastic localization, and simultaneous localization and mapping (SLAM)
are presented.
The algorithms are tested in numerical experiments where it is found that
fewer samples are required if implicit sampling is used, 
and that the overall runtimes of the algorithms are reduced.

\section{Introduction}
Many problems in physics and engineering require that one produces samples of a random variable with a given probability density function (pdf) $p$  \cite{ChorinHald,Doucet2001}. If $p$ is difficult to sample, one can sample an easier-to-sample pdf $p_0$ (the importance function), and attach to each sample $X^j \sim p_0$, $j=1,\dots,M$, the weight $W^j\propto p(X^j)/p_0(X^j)$ (capital letters denote realizations of random variables, and these realizations are 
indexed in superscript; subscript indices label discrete time). The weighted samples $X^j$ form an empirical estimate of the pdf $p$ and, under mild assumptions, the empirical estimate converges weakly to $p$, i.e., expected values of a function $h$ with respect to $p$ can be approximated by weighted averages of the function values $h(X^j)$. The difficulty is to find a ``good'' importance function~$p_0$: if $p_0$ is not a good approximation of $p$, for example, if $p_0$ is not large when $p$ is large, then many of the samples one proposes using $p_0$ may have a small probability with respect to $p$ and, therefore, carry little information, and make sampling inefficient. 

Implicit sampling is a general method for constructing useful importance functions \cite{chorintupnas,Chorin2010b}. 
The basic idea is to first locate the region of high probability with respect to $p$ via numerical optimization, and then to map a reference variable to this region; this mapping is done by solving algebraic equations. While the cost of generating a sample with implicit sampling is larger than in many other Monte Carlo (MC) schemes, implicit sampling can be efficient because the samples carry more information so that fewer samples are required. In \cite{atkins,Morzfeld2012,Morzfeld2011} this trade-off was examined in the context of geophysics, where it was found that implicit sampling indeed can be efficient, in particular when the dimension of the problem is large. 

Here implicit sampling is applied to estimation and control in robotics,
and new algorithms for stochastic optimal control, Monte Carlo localization, 
and simultaneous localization and mapping (SLAM) are developed.
Implementations and efficiencies of these algorithms are illustrated and explored with examples.
In particular, it is investigated if implicit sampling, which requires fewer samples that are however  
more informative and more expensive, can be efficient compared to other sampling schemes that 
may require more samples, each of which is less informative, however cheaper to generate.

The optimal control of a stochastic control problem can be found by solving a stochastic Hamilton-Jacobi-Bellman (HJB) equation \cite{Stengel}. The dimension of the domain of this partial differential equation (PDE) equals the dimension of the state space of the control problem. Classical PDE solvers require a grid on the domain and, therefore, are impractical for control problems of moderate or large dimension. For a class of stochastic optimal control problems, one can use MC solvers instead of grid based PDE techniques and, since the MC approach avoids grids, it is in principle feasible to solve larger dimensional control problems within this class \cite{Kappen2005a,Kappen2005b,Kappen2006,Theodorou2010}. However, the sampling scheme must be chosen carefully or else MC based PDE solving will also fail (in the sense that unfeasibly many samples are required). 
It is shown in section 3 how to apply implicit sampling to obtain a practical algorithm for stochastic optimal control
that avoids many of the pitfalls one faces in MC based PDE solving.
The method and algorithm are illustrated with the double slit problem (see section \ref{sec:DoubleSlit} and \cite{Kappen2005a}), which is a simple but vivid example of how things can go wrong, and how they can be fixed with implicit sampling. 
Using an inverse dynamics controller for a two-degrees-of-freedom robotic arm as an example,
it is also shown how to obtain a semi-analytic solution for a linear Gaussian problem via implicit sampling.
Finally, it is indicated how implicit sampling and stochastic optimal control can help with being trapped 
in local minima of non-convex optimization problems. 
An extension of the method presented here is also discussed in the conference paper \cite{Yang2014}.

In robotics, one often updates the predictions of a dynamic model of an autonomous robot with the output of the robot's sensors (e.g.~radar or lidar scans). This problem is often called ``localization'', and can be formulated as a sampling problem. Localization algorithms that rely on MC sampling for the computations are called Monte Carlo localization (MCL)~\cite{ThrunMCL}. One can also learn the geometry and configuration of the map while simultaneously localizing the robot in it, which leads to the problem of ``simultaneous localization and mapping''. Efficient numerical solutions of the MCL and SLAM problems are a fundamental requirement for autonomous robots \cite{FastSLAM,ThrunBook}. The various MCL and SLAM algorithms differ in their importance function $p_0$, and some algorithms have been shown to be inefficient due to a poor choice of $p_0$ \cite{FastSLAM}. Here it is shown how to use implicit sampling 
to generate an efficient importance function for MCL and SLAM.
The implicit sampling based MCL and SLAM algorithms are convergent, i.e.~as the number of samples goes to infinity
one obtains an empirical estimate of the true posterior even if the underlying dynamics are nonlinear
(whereas many other SLAM algorithms require linearity for convergence \cite{FastSLAM}).
The memory requirements of the implicit sampling based SLAM algorithm scale linearly 
with the number of features in the map.
The efficiencies of the new MCL and SLAM algorithms are compared 
to the efficiencies of competing MCL and SLAM algorithms in numerical experiments with the data set \cite{UTE_DataSet}.

\section{Review of implicit sampling}\label{sec:Review}
The efficiency of MC sampling depends on how well the importance function $p_0$ approximates the target pdf~$p$. In particular, $p_0$ must be large where $p$ is large, or the samples one produces using $p_0$ are unlikely with respect to $p$. Implicit sampling is a general MC sampling technique that constructs an importance function that is large where $p$ is large by mapping a reference variable to the region where $p$ is large
\cite{chorintupnas,Chorin2010b,Morzfeld2011,Morzfeld2012,atkins}. Here implicit sampling is described in general terms; below implementations of implicit sampling are described in the context of applications.

The region where $p$ is large is the region where $-\log(p)$ is small (the logarithm is used here because in applications $p$ often involves exponential functions). The region where $p$ is large can thus be located via minimization of the function
\begin{equation*}
	F(x)=-\log\left(p\left(x\right)\right).
\end{equation*}
If $F$ is convex, the minimizer $\mu = \mbox{argmin } F$ (i.e.~the mode of $p$) is the most likely value, and high-probability samples are in its neighborhood. One can obtain samples in this region by solving the stochastic algebraic equation
\begin{equation}
\label{eq:ImplicitSamplingEquation}
	F(x)-\phi = \frac{1}{2}\xi^T\xi,
\end{equation}
where $\phi = \min F$, $\xi\sim\mathcal{N}(0,I_m)$ is an $m$-dimensional Gaussian reference variable and where $T$ denotes a transpose and $I_m$ is the identity matrix of order~$m$; here, and for the remainder of this paper, $\mathcal{N}(a,B)$ is shorthand notation for  a Gaussian pdf with mean $a$ and covariance matrix~$B$. Note that the right-hand-side of (\ref{eq:ImplicitSamplingEquation}) is likely to be small because $\xi$ is close to the origin, which implies that the left-hand-side is small, which in turn means that the solution of (\ref{eq:ImplicitSamplingEquation}), i.e.~the sample, is close to the mode and, thus, in the high-probability region.

The weights of the samples are proportional to the absolute value of the Jacobian determinant
of the map that connects the sample $X^j$ to the reference variable $\xi$:
\begin{equation*}
	W(X^j) \propto \left|\det\left(\frac{\partial x}{\partial \xi}(X^j)\right)\right|.
\end{equation*} 
In practice, the weights are usually normalized so that their sum is one.
Various ways of constructing a map from $\xi$ to $x$ have been reported in the literature~\cite{Chorin2010b,Morzfeld2011}
two of which are summarized below.
In summary, implicit sampling requires (\emph{i}) minimizing $F=-\log p(x)$; and (\emph{ii}) solving equation (\ref{eq:ImplicitSamplingEquation}). This two-step approach makes efficient use of the available computational resources: the minimization identifies the high probability region and the samples are focused to lie in this region, so that (almost) all samples carry information. This is in contrast to other sampling schemes where an importance function is chosen ad-hoc, which often means that many of the samples carry little or no information; the computations used for generating uninformative samples are wasted. 

Finally, the assumption that $F$ is convex can be relaxed. Implicit sampling can be used without modification if $F$ is merely $U$-shaped, i.e.~the target pdf, $p$, is unimodal. Multi-modal target pdfs can be sampled e.g.~by using mixture models, for which one approximates each mode using the above recipe (see \cite{Chorin2010b,Morzfeld2011,Brad2013} for more detail).

\subsection{Generating samples with random maps}
To generate samples, one solves (\ref{eq:ImplicitSamplingEquation}) with a one-to-one and onto mapping.
There are many choices for such a mapping, and one is to solve (\ref{eq:ImplicitSamplingEquation}) in 
a random direction, i.e.~one puts
\begin{equation}
	\label{eq:RandomMap}
	X^j = \mu+\lambda^j L^{-T}\eta^j,
\end{equation}
where $L$ is a Cholesky factor of the Hessian of $F$ at the minimum, $\eta^j = \xi^j/\sqrt{(\xi^j)^T\xi^j}$ is a vector which is uniformly distributed on the $m$-dimensional hypersphere, and where $\lambda^j$ is a scalar.
One then computes $\lambda^j$ by substituting~(\ref{eq:RandomMap}) into~(\ref{eq:ImplicitSamplingEquation}). This approach, called ``random map'', requires solving a scalar equation to generate a sample~\cite{Morzfeld2011}. Moreover, the Jacobian of this map is easy to evaluate with the formula (see \cite{Morzfeld2011} for a derivation)
\begin{equation}
\label{eq:Jacobian}
	\left|\det\left(\frac{\partial x}{\partial \xi}\right)\right|= \frac{(\xi^T\xi)^{(1-m)/2}}{\det\left(L\right)} \left|\frac{\lambda(\xi)^{m-1}}{2\nabla F\; L^{-T}\eta}\right|,
\end{equation}
where $\nabla F$ is the gradient of $F$ with respect to $x$ (an $m$-dimensional row vector).

In numerical implementations of this method, calculating $\nabla F$ may require repeated evaluations of $F$, e.g.~if the gradient is approximated via
finite differences. 
This can be avoided by noting that
\begin{equation*}
	\frac{d\lambda}{d\rho} = \frac{1}{2\nabla F\; L^{-T}\eta},
\end{equation*}
where $\rho=\xi^T\xi$, so that the Jacobian becomes
\begin{equation}
\label{eq:Jacobian2}
	\left|\det\left(\frac{\partial x}{\partial \xi}\right)\right|= \frac{\rho^{(1-m)/2}}{\det\left(L\right)} \left|\lambda(\xi)^{m-1}\frac{d\lambda}{d\rho}\right|.
\end{equation}
The scalar derivative $d\lambda/d\rho$ can be evaluated using finite differences
with a few evaluations of $F$ (the precise number of function evaluations
depends on the finite difference scheme one chooses).

\subsection{Generating samples with approximate quadratic equations}
Another path to solving (\ref{eq:ImplicitSamplingEquation}) is to replace it 
with an approximate quadratic equation 
\begin{equation}
\label{eq:AppSampleEqn}
	\hat{F}(x)-\phi = \frac{1}{2}\xi^T\xi,
\end{equation}
where
\begin{equation*}
\hat{F}(x)=\phi+\frac{1}{2}\left(x-\mu\right)^TH\left(x-\mu\right),
\end{equation*} 
is the Taylor expansion of order two of $F$ about its minimizer~$\mu$, and where $H$ is the Hessian of $F$ at the minimum. A solution of (\ref{eq:AppSampleEqn}) is
\begin{equation*}
	X^j = \mu+L^{-T}\xi_j,
\end{equation*}
where $L$ is the Cholesky factor of $H=LL^T$. The weights are
 \begin{equation}
 \label{eq:WeightsQuadratic}
	W^j\propto  \exp\bigl(-\phi+\hat{F}(X^j)-F(X^j)\bigl),
\end{equation}
and account for the error one makes by solving (\ref{eq:AppSampleEqn}) instead of (\ref{eq:ImplicitSamplingEquation})  (see \cite{Chorin2010b} for more detail).

\section{Application to path integral control}
The finite horizon stochastic optimal control problem considered here is as follows: find a control $u$ (a $p$ dimensional vector function of the state $x$) that minimizes the cost function
\begin{equation}
	C(x_0,t_0,u) = E\left[\Phi(x(t_f)) + \displaystyle\int_{t_0}^{t_f}u(x(\tau),\tau)^TRu(x(\tau),\tau)+V(x(\tau),\tau)d\tau\right],\nonumber
\end{equation}
where $x$ is an $m$-dimensional vector (the state), $t_0\leq t\leq t_f$ is time, $t_f$ is the final time (or the horizon), $\Phi$ is a scalar function that describes the ``final cost'', $R$ is a $p\times p$ symmetric positive definite matrix that specifies the control cost, and $V$ is a scalar function which describes the ``running cost'' (which is also called ``potential''); the expectation is taken over trajectories of the stochastic differential equation (SDE) 
\begin{equation}
\label{eq:StochControlDyanamics}
	dx = f(x,t)\,dt+G\left(u\,dt+Q\,dW\right),
\end{equation}
starting at $x(t_0)=x_0$, where $f$ is a smooth $m$-dimensional vector function which describes the dynamics, and where $G$ and $Q$ are $m\times p$ and $p\times r$ matrices describe how the uncertainty and controls 
are distributed within the system; $W$ is Brownian motion (see, e.g.~\cite{Stengel} for more detail about stochastic optimal control; I closely follow \cite{Kappen2005a}
in the presentation of path integral control). 

The ``optimal cost-to-go'' is defined as
\begin{equation}
  J(x,t) = \displaystyle\min_u C(x,t,u), 
  \nonumber
\end{equation}
and satisfies the stochastic Hamilton-Jacobi-Bellman (HJB) equation \cite{Stengel}:
\begin{eqnarray}
  \partial_t J(x,t)& =& \partial_xJ(x,t) +\mbox{Tr}\left(Q^TG^T\partial_{xx}J(x,t)GQ\right)+V(x,t) \nonumber \\ 
  & &-\tfrac{1}{2}\partial_xJ(x,t)G^TR^{-1}G\partial_xJ(x,t),\nonumber
\end{eqnarray}
where $\mbox{Tr}$ is the trace. If there exists a scalar $\gamma$ such that 
\begin{equation}
\label{eq:NoiseCostRequirement}
  \gamma GR^{-1}G^T=GQQ^TG^T,
\end{equation}
then the nonlinear change of variables 
\begin{equation}
  J(x,t)=-\gamma\log \psi(x,t),
  \nonumber
\end{equation}
linearizes the stochastic HJB equation and one obtains
\begin{equation}
  \label{eq:HJBForward}
  \partial_t\psi=\frac{V(x,t)}{\gamma}\psi-f(x,t)^T\partial_x\psi-\tfrac{1}{2}\mbox{Tr}\left(Q^TG^T\partial_{xx}\psi GQ\right),
\end{equation}
with final condition $\psi(x,t_f)=\exp(-\Phi(x(t_f))/\gamma)$ and with the optimal control
\begin{equation}
\label{eq:optimalControl}
  u=-R^{-1}\partial_x J \, G,
\end{equation}
see \cite{Fleming,Kappen2005a}. Thus, one can compute the optimal control by solving the HJB equation. Numerical PDE solvers typically require a grid on the domain of the PDE, however the domain has the dimension of the state space of the control problem. Thus, grid based numerical PDE techniques only apply to control problems of a low dimension
(or else the computational requirements become excessive).

For larger dimensional problems, one can use stochastic PDE solvers which do not require a grid. In particular, one  can solve the adjoint equation
\begin{equation}
  \partial_t\psi=-\frac{V(x,t)}{\gamma}\psi-\partial_x(f(x,t)\psi)+\tfrac{1}{2}\mbox{Tr}\left(QG\partial_{xx}\psi G^TQ^T\right),
  \nonumber
\end{equation}
forward in time (instead of solving (\ref{eq:HJBForward}) backwards in time) with the Feynman-Kac formula (see, e.g.\cite{ChorinHald}) 
\begin{equation}
  \label{eq:Integral}
  \psi(x,t)= E\left[ \exp\left(-\frac{1}{\gamma}\left(\Phi(y(t_f),t_f)+\int_t^{t_f}V(y(\tau),\tau)d\tau \right)\right)   \right],
\end{equation}
where the expectation is over the trajectories of
\begin{equation}
  \label{eq:PathIntegralDynamics}
  dy=f(y,\tau)\,d\tau+GQ\,dW,
\end{equation}
starting at $y(\tau=t)=x(t)$.
This is the path integral formulation of stochastic optimal control \cite{Kappen2005a,Kappen2005b,Kappen2006,Theodorou2010}. Note that the class of problems that can be tackled with path integral control is rather general, since the assumptions of (\emph{i}) a quadratic control cost; and (\emph{ii}) the condition on the noise and the control cost in equation (\ref{eq:NoiseCostRequirement}) are not restrictive. In particular, the dynamics $f(x,t)$ and the potential $V(x,t)$ in (\ref{eq:StochControlDyanamics}) can be nonlinear.

To evaluate the Feynman-Kac formula numerically, one approximates the infinite dimensional integral in (\ref{eq:Integral}) by a finite dimensional one. For example, one can discretize the integral on a regular grid in time with constant time step $\Delta t$, so that  
\begin{eqnarray}
\label{eq:DiscreteExpectation}
  \psi(x,t) &\approx& \int dy_1 \cdots \int dy_n  \; p(y_1,\dots,y_n) \nonumber \\
  & &\times\exp\left(-\frac{1}{\gamma}\Phi(y_n,n\Delta t) -\frac{\Delta t}{2\gamma}\sum_{i=1}^nV(y_i,\tau_i)+V(y_{i-1},\tau_{i-1})\right),  
\end{eqnarray}
where $y_1=x(t)$, $n=(t_f-t)/\Delta t$, $\tau_i=t+i\Delta t$, $i=0,1,\dots,n$, and where $p$ is the pdf of the discretized trajectory $y_1,\dots,y_n$. Here the trapezoidal rule is used to discretize the integral of the potential, however other choices are also possible~\cite{Hammond1994}. The SDE (\ref{eq:PathIntegralDynamics}) implies that the increments $\Delta y_i = y_i-y_{i-1}$ are independent Gaussians, e.g.~for a forward Euler discretization~\cite{Kloeden} of (\ref{eq:PathIntegralDynamics}),
\begin{equation}
  \Delta y_i \sim \mathcal{N}\left(f(y_i,\tau_i)\Delta t,\Sigma\Delta t\right),
  \nonumber
\end{equation}
where $\Sigma = \gamma GR^{-1}G^T=GQQ^TG^T$, so that
\begin{eqnarray}
\label{eq:ControlF}
  p(y_1,\dots,y_n) &=& \prod_{i=1}^n p(\Delta y_i)\nonumber\\
  & =&\frac{\exp\left(-\frac{1}{2}\sum_{i=1}^n\Delta t\left(\frac{\Delta y_i}{\Delta t}-f(y_i,\tau_i)\right)^T\Sigma^{-1}\left(\frac{\Delta y_i}{\Delta t}-f(y_i,\tau_i)\right)\right)}{\left(2\pi \Delta t\det \Sigma \right)^{-n/2}}.
  \end{eqnarray}
Thus,
\begin{equation}
  \psi(x,t) \approx \frac{1}{\left(2\pi \Delta t\det \Sigma \right)^{n/2}}\int dy_1 \cdots \int dy_n \exp\left(-F(y_1,\dots,y_n)\right),  
  \nonumber
\end{equation}
where
\begin{eqnarray}
\label{eq:ContrlFII}
F(y_1,\dots,y_n)&=&\frac{1}{\gamma}\Phi(y_n,n\Delta t) +\frac{\Delta t}{2\gamma}\sum_{i=1}^nV(y_i,\tau_i)+V(y_{i-1},\tau_{i-1})
  \nonumber \\
  & & +\frac{1}{2}\sum_{i=1}^n\Delta t\left(\frac{\Delta y_i}{\Delta t}-f(y_i,\tau_i)\right)^T\Sigma^{-1}\left(\frac{\Delta y_i}{\Delta t}-f(y_i,\tau_i)\right).
\end{eqnarray}
Note that this $F$ depends on how one discretizes the integral of the potential $V$, and the SDE (\ref{eq:PathIntegralDynamics}); other discretization schemes will lead to different~$F$s. 

For a given discretization, i.e.~for a given $F$, MC sampling can be applied to compute the expectation in (\ref{eq:DiscreteExpectation}). For example, one can evaluate 
\begin{equation*}
G(y_1,\dots,y_n) = \frac{1}{\gamma}\Phi(y_n,n\Delta t) +\frac{\Delta t}{2\gamma}\sum_{i=1}^nV(y_i,\tau_i)+V(y_{i-1},\tau_{i-1}),
\end{equation*} 
for discretized trajectories of (\ref{eq:PathIntegralDynamics}), followed by averaging. However, this method can fail if the potential $V$ has deep wells \cite{Kappen2005a,Kappen2005b,Kappen2006}. In this case, many trajectories of (\ref{eq:PathIntegralDynamics}) can end up where $V$ is large and, thus, contribute little to the approximation of the expected value $\psi$. For efficient sampling, one needs a method that guides the samples to remain where the potential is small. 

This guiding of the samples can be achieved via implicit sampling.
Recall that in implicit sampling one first locates the region of high probability by minimizing
$F=-\log (p(x))$, where $p(x)$ is the pdf of the random variable one is interested in.
The trajectories we wish to sample have the joint pdf (\ref{eq:ControlF}),
so that, in order to apply implicit sampling to path integral control, one needs to minimize $F$ in (\ref{eq:ContrlFII}).
With this $F$, one solves the algebraic equations~(\ref{eq:ImplicitSamplingEquation}) with a one-to-one and onto map (\ref{eq:RandomMap}). The integral~(\ref{eq:DiscreteExpectation}) becomes
\begin{eqnarray}
\label{eq:ImplicitSamplingSolution}
  \psi(x,t) &\approx& \frac{\exp(-\phi)}{\left(2\pi \Delta t\det \Sigma\right)^{n/2}} \int d\xi_1 \cdots \int d\xi_n \exp(-\xi^T\xi/2)\left|\frac{\partial x}{\partial \xi}\right|\nonumber\\
  &\approx&\frac{\exp(-\phi)}{\left(\Delta t\det \Sigma\right)^{n/2}}\; E_\xi\left[\left|\frac{\partial x}{\partial \xi}\right|\right],
\end{eqnarray}
where the expectation is over the reference variable $\xi$ and where $\partial x/\partial \xi$ is the Jacobian of the 
map from $x$ to $\xi$ in equation (\ref{eq:RandomMap}). 
Combining (\ref{eq:ImplicitSamplingSolution}) with the expression (\ref{eq:Jacobian}) for the Jacobian gives
\begin{equation*}
\psi(x,t)\approx \frac{\exp(-\phi)}{2\left(\Delta t\det \Sigma\right)^{n/2}\det\left(L\right)}\; 
E_\xi\left[(\xi^T\xi)^{(1-m)/2} \left|\frac{\lambda(\xi)^{m-1}}{\nabla F(\xi) \,\eta}\right|\right].
\end{equation*}
The expected value is now straightforward to compute via Monte Carlo,
i.e.~sampling a the reference variable $\xi$ to obtain $M$ samples $\xi^j$, $j=1,\dots,M$, computing, for each one,
$\lambda(\xi^j)$ and the gradient of $F$, followed by averaging.
In numerical implementations, it may be more efficient to use the equivalent expression
(\ref{eq:Jacobian}) for the Jacobian (which avoids computing the gradient of $F$).
In this case, one obtains
\begin{equation*}
\psi(x(t),t)\approx \frac{\exp(-\phi)}{\left(\Delta t\det \Sigma\right)^{n/2}\det\left(L\right)}\; 
E_\xi\left[(\xi^T\xi)^{(1-m)/2} \left|\lambda(\xi)^{m-1}\frac{d\lambda}{d\rho}(\xi)\right|\right].
\end{equation*}

Once $\psi(x,t)$ is computed, 
we can compute the optimal control from (\ref{eq:optimalControl}), e.g.~via finite differences.
It is important to note that this defines the optimal control at time $t$,
and that the computations have to be repeated at the next time step to compute $\psi(x(t+dt),t+dt)$.
The use of implicit sampling in stochastic optimal control is illustrated with three examples.

\subsection{The double slit problem}\label{sec:DoubleSlit}
The double slit problem is a simple example that illustrates the pitfalls one must avoid
when using MC sampling for path integral control \cite{Kappen2005a}.
The problem is as follows. Suppose you observe a somehow confused person walking (randomly) towards a wall with two doors, and your job is to guide this person through either one of the two doors. The person and you are modeled by the controlled dynamics
\begin{equation}
  dx=u\,dt+\sqrt{\sigma}\,dW, \nonumber
\end{equation}
the final cost is quadratic, $\Phi = x(t_f)^2/2$, and the scalar $R>0$, 
that defines the cost of the control, is given; the ``wall'' is modeled by the potential
\begin{equation}
V(x,t) = \left\{
        \begin{array}{ll}
            \infty & \mbox{if } t=t_1 \mbox{ and }x<a, \mbox{ or } b<x<c, \mbox{ or } d<x  \\
            0 & \mbox{otherwise},
        \end{array}
    \right. \nonumber
\end{equation}
for given $a<b<c<d$ and $t_1>0$. 
The stochastic HJB equation of the double slit problem is
\begin{equation}
	\psi_t=-\frac{Q}{\gamma}+\frac{1}{2}\sigma \psi_{xx}, \nonumber
\end{equation}
where $\gamma=r\sqrt{\sigma}$ satisfies (\ref{eq:NoiseCostRequirement}). 

\begin{figure}[tbp]
\begin{center}
{\includegraphics[width=1\textwidth]{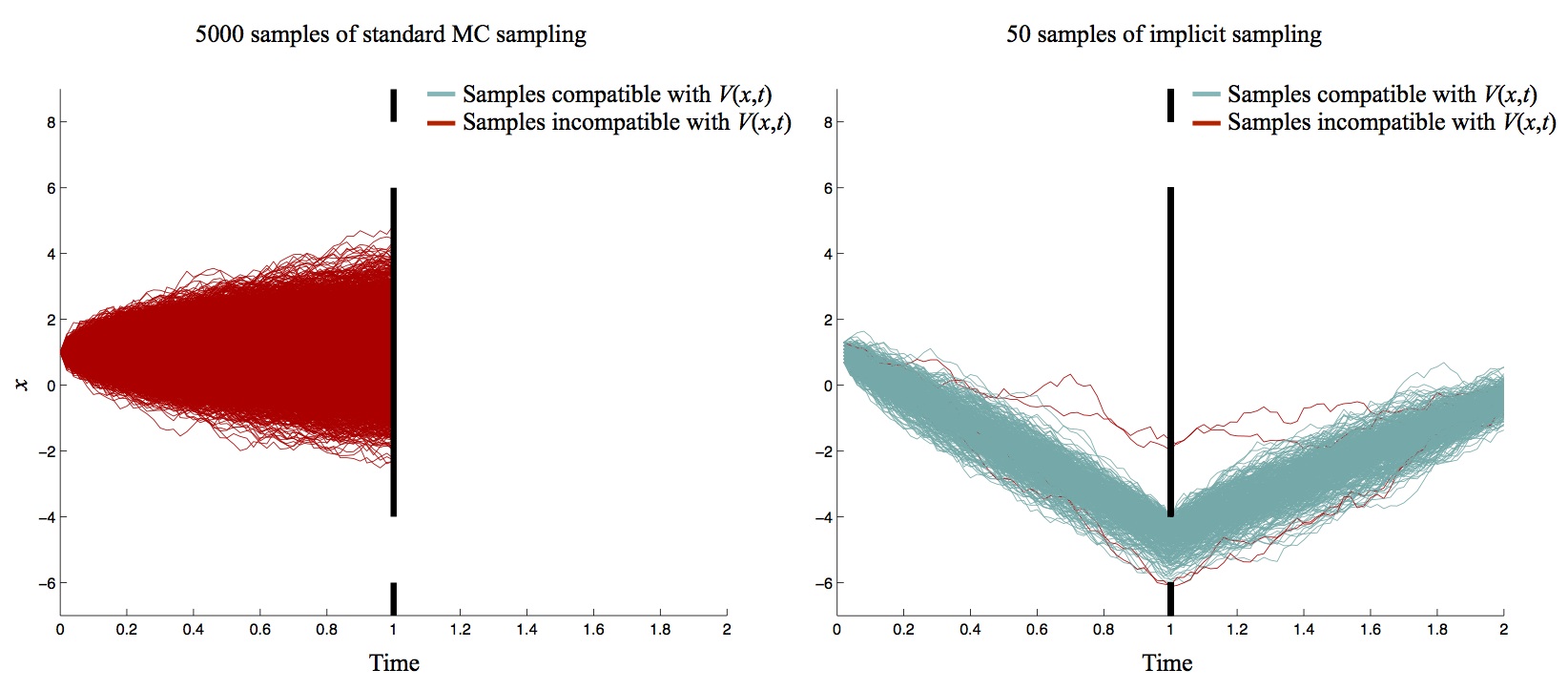}}
\caption{Random walks for the double slit problem. Left: standard MC sampling with 5000 unguided random walks, all of which hit the potential wall. Right: 50 guided random walks obtained via implicit sampling, only 4 of which hit the potential wall.} 
\label{fig:DoubleSlitWalks} 
\end{center}
\end{figure}

One can attempt to solve this equation using standard MC sampling as follows. Solve the SDE
\begin{equation}
	\label{eq:DoubleSlitSDE}
	dy=\sqrt{\sigma}\,dW, \nonumber
\end{equation}
starting at $y(\tau=t)=x(t)$ repeatedly, e.g.~using a uniform grid in time and the forward Euler scheme~\cite{Kloeden}
\begin{equation}
	y_{n+1}=y_{n}+\sqrt{\sigma \Delta t}\Delta w^n, \nonumber
\end{equation}
where $\Delta w^n$ are independent Gaussians with mean 0 and variance 1. For each trajectory $\left\{y_1,\dots,y_n\right\}$, evaluate 
\begin{equation}
 G(y_1,\dots,y_n)=\frac{1}{2\gamma}y_n^2 +\frac{\Delta t}{2\gamma}\sum_{i=1}^nV(y_i,\tau_i)+V(y_{i-1},\tau_{i-1}).
  \nonumber
\end{equation}
With a high probability, the trajectories will hit the potential wall at $t_1$ and, thus, $G$ is infinite, so that the contribution of a trajectory to the expected value $\psi(x,t)$ in(\ref{eq:Integral}) is zero with a high probability. The situation is illustrated with a simulation using the parameters of table~\ref{tab:DoubleSlitParameters}. 
\begin{table}[tb]
\caption{Parameters of the double slit problem}
\begin{center}
\begin{tabular}{l c c}
Description & Parameter & Value\\
\hline
Final time &$t_f$ & $2$\\ 
Critical time &$t_1$ & $1$\\ 
Slits & $a,b,c,d$ & $-6,-4,6,8$\\
Initial position &$x(0)$ & 1\\
Variance of the noise &$\sigma$&$1$\\
Control cost& $r$&$0.1$\\
Time step & $\Delta t$ & 0.02\\
\hline
\end{tabular}
\end{center}
\label{tab:DoubleSlitParameters}
\end{table}
The results are shown in the left panel of figure~\ref{fig:DoubleSlitWalks}.
One observes that 5000 unguided random walks from (\ref{eq:DoubleSlitSDE}), all starting at $x(0)$, hit the potential wall, and, thus, score $G=\infty$. All 5000 samples thus carry no probability and do not contribute to the approximation of $\psi$. Even when $O(10^5)$ samples are considered, it is unlikely that sufficiently many make it past the potential wall\cite{Kappen2005a}. We conclude that this method is unfeasible for this problem, since the number of samples required is extremely large due to the deep wells in the potential.

Implicit sampling can be applied here to fix these problems. In particular, one finds from the Feynman-Kac formula that
\begin{equation}
	\psi(x,t) = E\left[\exp\left(-\frac{y(t_f)^2}{2\gamma}-\frac{1}{\gamma}\int_t^{t_f}V(y(\tau),\tau)d\tau\right)\right],\nonumber
\end{equation}
which one can compute with implicit sampling in two parts. For $t\geq t_1$, i.e.~after one has passed the potential wall, the probability of each path is Gaussian with variance $s=t_f-t$, so that  
\begin{equation}
	\psi(x,t) = \int \frac{1}{\sqrt{2\pi s}}\exp\left(-\frac{(z-y)^2}{2s}\right) \exp\left(-\frac{z^2}{2\gamma}\right)dz \quad \mbox{for } t\geq t_1.\nonumber
\end{equation}
For implicit sampling, one defines
\begin{equation}
	F(y) = \frac{(z-y)^2}{2s}+\frac{z^2}{2\gamma},\nonumber
\end{equation}
whose minimizer and minimum are $\mu =\mbox{argmin }F = y \gamma/(s+\gamma)$, $\phi = \min F = y^2/(s+\gamma)$. Since $F$ is quadratic, $F = \phi +H(y-\mu^2)/2$, where $H =(s+\gamma)/(s\gamma)$ is the Hessian of $F$ at the minimum. The algebraic equation (\ref{eq:ImplicitSamplingEquation}) can thus be solved by the linear map
\begin{equation}
	y = \mu + \sqrt{\frac{s\gamma}{s+\gamma}}\xi,\nonumber
\end{equation}
so that 
\begin{equation}
	\psi(x,t) = \exp(-\phi)\sqrt{\frac{\gamma}{(s+\gamma)}}, \quad \mbox{for } t\geq t_1.\nonumber
\end{equation}

For $t<t_1$, one splits the integration into two parts, first going from time $t$ to $t_1$ with probability $p(y_1,t_1|y,t)$, and then from $t_1$ onwards to $t_f$ with probability $p(z,t_f|y_1) = \mathcal{N}(y_1,s)$, $s=t_f-t_1$:
\begin{equation}
	\psi(x,t) = \int dz\int dy_1 \exp\bigl(-\frac{z^2}{2\gamma}\bigl)p(z,t_f|y_1)p(y_1,t_1|y,t).\nonumber
\end{equation}
Implicit sampling uses the information about the potential which is infinite at $t_1$ except for the two slits, so that $p(y_1,t_1|y,t)=0$ outside the slits and Gaussian with mean $y$ and variance $\hat{s}=t_1-t$ otherwise. Since only the slits need to be explored with samples, the integration can be carried out over the slits:
\begin{equation}
	\psi(x,t) = \int dy\displaystyle \int_{y\in [a,b]\times[c,d]} dy_1 \exp\bigl(-\frac{z^2}{2\gamma}\bigl)\exp\bigl(-\frac{(z-y_1)^2}{2s}\bigl)\exp\bigl(-\frac{(y_1-y)^2}{2\hat{s}}\bigl).\nonumber
\end{equation}
The evaluation of the above integral using the same strategy as above for $t\geq t_1$ is tedious, but straightforward. Note that the implicit sampling strategy here is the key to solving this problem, because it locates the wells in the potential.

The above analytic solution is compared to a numerical implementation of implicit sampling, for which the paths are discretized with a constant time step $\Delta t$. Here the numerical integration is also split up into two parts. For $t\geq t_1$ the function $F$ is quadratic because the probabilities are Gaussian and the potential has no effect. Thus,
\begin{equation}
	F = \frac{y_n^2}{2\gamma}+\frac{1}{2\gamma \Delta t}\sum_{i=1}^n(y_i-y_{i-1})^2,\nonumber
\end{equation}
where $n=(t_f-t)/\Delta t$ and $y_i=y(t+i\Delta t)$,$i=1,\dots,n$. Minimizing $F$ is straight forward, and the algebraic equation (\ref{eq:ImplicitSamplingEquation}) can be solved with a linear map
\begin{equation}
	y = \mu+L^{-T}\xi,\nonumber
\end{equation}
where $y=(y_0,\dots,y_n)$ is a $n$-dimensional vector, $\mu$ is the minimizer of $F$ and $\xi$ is an $n$-dimensional vector whose elements are independent standard normal variates; $L$ is a Cholesky factor of the Hessian at the minimum. Since the Jacobian of this linear map is constant, (\ref{eq:ImplicitSamplingSolution}) gives
\begin{equation}
	\psi(x,t) = \frac{\exp(-\phi)}{(\gamma \Delta t)^{n/2}\det L},\nonumber
\end{equation}
where $\phi$ is the minimum of $F$. 

For $t<t_1$, one obtains the same $F$, but needs to perform a constraint minimization over the slits. There are two (local) minima, one per slit, and both can be found easily using quadratic programming \cite{Nocedal,Fletcher}. One can then generate a sample close to each of the minima using 
\begin{equation}
	y= \mu+L^{-T}\xi_j,\nonumber
\end{equation}
where $y$ again is a vector whose elements are the discretized path and where $\mu$ is the location of a local minimum of the constraint problem and $L$ is a Cholesky factor of the unconstrained Hessian at a minimum. 
Equation~(\ref{eq:ImplicitSamplingSolution}) becomes
\begin{equation}
\psi(x,t) \approx \frac{\exp(-\phi)}{\left(\Delta t\gamma\right)^{n/2}} \int_{\mbox{Slits}} \frac{1}{\det L}d\xi.\nonumber
\end{equation}
The expected value of the Jacobian ($1/\det L$) over the slits can be computed by Monte Carlo as
follows. Generate $M$ samples $\xi^j$, $j=1,\dots,M$, and, for each one, compute the
corresponding trajectory and check if it hits the potential wall. 
The integral is $1/\det(L)$ times the ratio of the number of trajectories that pass through the wall
and $M$. 

The right panel of figure \ref{fig:DoubleSlitWalks} shows 50 trajectories one obtains with this approximation at $t=\Delta t$. Note that most of the trajectories pass through the slit, i.e.~most of the samples carry a significant probability, score a small $F$, and, thus, contribute to the approximation of the expected value $\psi(x,t)$. These ``guiding'' effects make it possible to solve the problem with $O(10)$ samples, while the standard Monte Carlo scheme fails even with $O(10^5)$ samples \cite{Kappen2005a}. The Laplace guided strategy presented in \cite{Kappen2005a} (which uses similar ingredients as implicit sampling and is also related to the optimal nudging constructions in \cite{Weare2012, Weare2013}) requires about $100$ samples.

\begin{figure}[tb]
\begin{center}
{\includegraphics[width=1\textwidth]{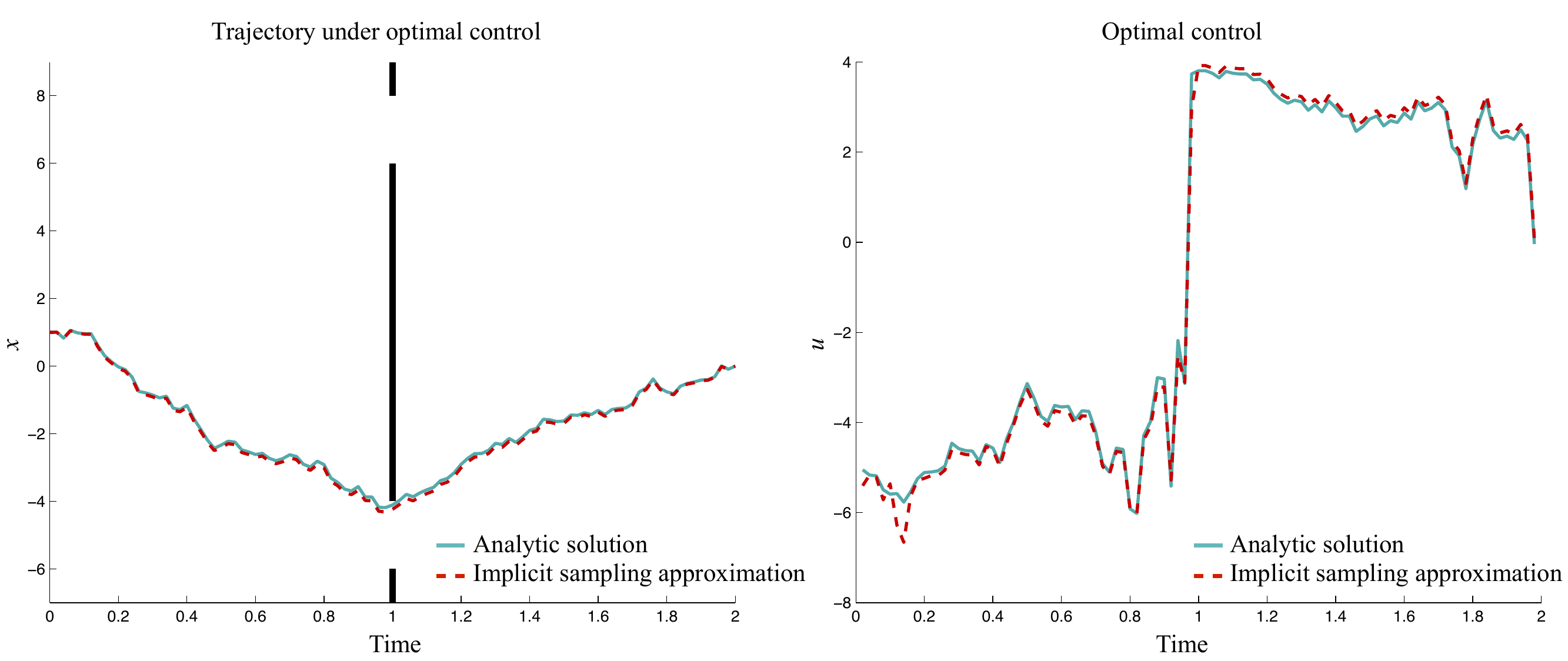}}
\caption{Comparison of the analytical solution and its approximation by implicit sampling with 50 samples. Left: the trajectory under optimal control (solid turquoise) and its numerical approximation (dashed red). Right: the optimal control (solid turquoise) and its numerical approximation (dashed red).} 
\label{fig:DoubleSlitControl} 
\end{center}
\end{figure}

The numerical approximation obtained with implicit sampling (50 samples) is compared to the analytical solution in figure \ref{fig:DoubleSlitControl}, where the optimal control and the trajectory under this control are shown.
One observes that the numerical approximation of the optimal control is quite close to the control one obtains from the analytical solution (right panel of figure \ref{fig:DoubleSlitControl}) and, hence, the controlled trajectory one obtains with implicit sampling is also close to the one computed from the analytical solution (left panel of figure \ref{fig:DoubleSlitControl}). This statement can be made more precise by solving the control problem repeatedly (each solution is a random event) and averaging. 

The double slit problem is solved 500 times and the error of the numerical approximation is recorded for each run. The Euclidean norm of the difference of the analytical solution and the approximation via implicit sampling is used to measure an error, in particular in $x$ (the trajectory under optimal control) and $u$ (the optimal control). The number of samples of implicit sampling is varied to study the convergence of the algorithm. The results are shown in table \ref{tab:DoubleSlitResults} where the mean and standard deviation of the Euclidean norm of the error in $x$, respectively $u$, scaled by the mean of the norm of $x$ and $u$ respectively, are listed.
\begin{table}[tb]
\caption{Mean errors of the numerical implementation of implicit sampling}
\begin{center}
\begin{tabular}{c c c c c}
 & \multicolumn{3}{c}{Number of samples} \\
 \hline
 & 10 & 50 & 100 & 500\\
Error in $x$ in \%& 2.64$\pm$0.63 & 2.70$\pm$0.63 &2.66$\pm$0.63 &2.61$\pm$0.62\\
(mean$\pm$std.~dev.) & & & &  \\
Error in $u$ in \%& 4.81$\pm$2.15 & 4.61$\pm$1.89 &4.54$\pm$2.37 &4.46$\pm$2.17\\ 
(mean$\pm$std.~dev.) & & & &  \\
\hline
\end{tabular}
\end{center}
\label{tab:DoubleSlitResults}
\end{table}
These numbers indicate the errors one should expect in a typical run. The errors are relatively small when 50 samples are used. The errors do not dramatically decrease for larger sample numbers, which indicates that the algorithm has converged. The small variance of the errors indicates that similar errors occur in each run, so that one concludes that implicit sampling is reliable.

The error in the numerical implementation of implicit sampling is mostly due to neglecting one of the local minima of $F$. In numerical tests, only a slight improvement was observed when both minima were used for sampling, however the computations are twice as fast if one considers only the smaller of the two minima.

\subsection{Stochastic control of a two-degrees-of-freedom robotic arm}
\begin{figure}[tb]
\begin{center}
{\includegraphics[width=0.25\textwidth]{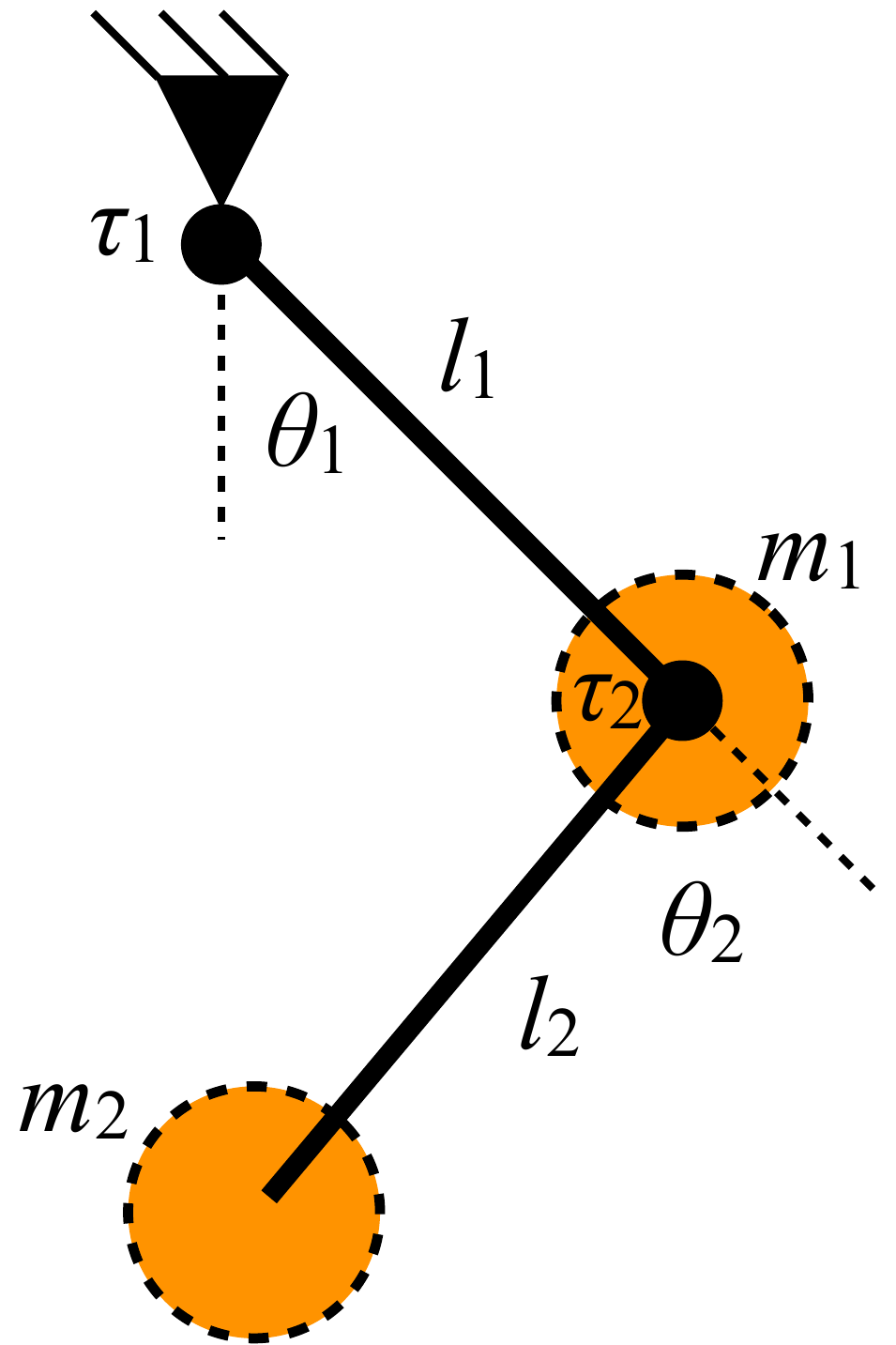}}
\caption{Sketch of a two-degrees-of-freedom robotic arm.} 
\label{fig:RobotArmSketch} 
\end{center}
\end{figure}

Stochastic optimal control of the two-degrees-of-freedom robotic arm shown in figure~\ref{fig:RobotArmSketch}
is considered. 
The controller is an ``inverse dynamics controller'' that linearizes the system.
This example demonstrates how implicit sampling based
path integral control works in linear problems.
Specifically, a semi-analytical solution is derived for which a numerical 
optimization is required, however the expected value
of the implicit sampling solution in (\ref{eq:ImplicitSamplingSolution})
can be evaluated explicitly (i.e.~no numerical sampling is needed).
Moreover, only some of the dynamic equations 
of the first order formulation of the dynamics are driven by noise,
so that the dimension of the stochastic subspace is less than the 
actual dimension of the problem.
Implicit sampling can exploit this features, and this example
demonstrates how.
The algorithm is tested in numerical experiments in which false
parameters are given to the algorithm to demonstrate the robustness 
of the path integral controller to model error.

The goal is to compute two independent torques $\tau_1$ and $\tau_2$ that drive the arm to a desired position $\theta^*$ (described by the two-angles $\theta_1,\theta_2$, the degrees-of-freedom). Neglecting energy dissipation, and assuming the robot is mounted on a horizontal table (no effects of gravity), the equation of motion is
\begin{equation}
	M(\theta)\ddot{\theta}+C(\theta,\dot{\theta})\dot{\theta} = \tau
	\nonumber
\end{equation}
where $\theta = (\theta_1,\theta_2)$ is a 2-dimensional vector and where dots denote derivatives with respect to time; the $2\times 2$ matrices \begin{eqnarray}
	M(\theta)&=& \left(\begin{array}{cc} 
	      a_1+a_3 \cos(\theta_2) & a_2+a_3\cos(\theta_2)  \\       
	      a_2+a_3\cos(\theta_2)  & a_2   	
	      \end{array}\right),\nonumber\\
	C(\theta,\dot\theta)&=&\left(\begin{array}{cc}
							-a_3 \sin(\theta_2)\dot{\theta}_2 & -a_3(\dot{\theta}_1+\dot{\theta}_2)\sin(\theta_2) \\
							 a_3 \sin(\theta_2)\dot{\theta}_1 & 0 \end{array}\right)\nonumber,
	\nonumber
\end{eqnarray}
where
\begin{eqnarray}
	a_1 & = & m_1 l^2_{c1} + m_2 l^2_{c1} + m_2 l^2_{c2} + I_1 + I_2\nonumber\\
	a_2 & = & m_2 l^2_{c2}+I_2\nonumber\\
	a_3 & = & l_1 m_2 l_{c2}\nonumber
\end{eqnarray}
define the dynamics and depend on the lengths of the arms $l_1,l_2$ and the loads $m_1, m_2$ (see, e.g.~\cite{Slotine1991} for more detail about this model and table \ref{tab:RoboticArm} for the parameters used in simulations). 
\begin{table}[tb]
\caption{Parameters of the robotic arm}
	\centering
	\begin{tabular}{ l c r}
    Description & Parameter & Value   \\ \hline
    Length of link one												 & $l_1$    & 1 m   		\\ 
    Distance to the center of mass of link one & $l_{c1}$ & 0.5 m  	\\
    Distance to the center of mass of link two & $l_{c2}$ & 0.5 m   \\
    Mass of link one													 & $m_1$ 		& 1 kg 		\\
    Mass of link two 													 & $m_2$	  & 1 kg			\\
    Moment of inertia of link one 						 & $I_1$ 		& 2 kg m$^2$	\\
    Moment of inertia of link two						 & $I_2$	  & 2 kg m$^2$ 	\\ \hline
    \end{tabular}
	\label{tab:RoboticArm}
\end{table}
One can compute the torques by inverting the dynamics and choosing $\tau = C(\theta,\dot\theta)+M(\theta)u$ to derive the linear system
\begin{equation}
\label{eq:RobotProblem}
	\left(\begin{array}{l}
            \dot\theta\\\ddot\theta
        \end{array}\right)=\left(\begin{array}{cc}
            0 & I\\
            0 & 0
        \end{array}\right)\left(\begin{array}{c}
            \theta\\
            \dot\theta
        \end{array}\right)+\left(\begin{array}{c}
            0\\
            I
        \end{array}\right)u,
\end{equation}
where $0$ denotes the matrix whose elements are all zero (of appropriate dimensions), and $I$ is the $2\times 2$ identity matrix; $u$ is a $2$-dimensional vector of controls which must be computed. To make the control robust to modeling errors, one can add noise and solve the stochastic optimal control problem
\begin{equation}
	dx=A\,dt+Gu\,dt+GQ\,dW,
	\nonumber
\end{equation}
where $x=(\theta,\dot\theta)$ is a 4-dimensional vector and where the matrices $A,G$ and $Q$ can be read from (\ref{eq:RobotProblem}), and where $W$ is a 2-dimensional Brownian motion. The final cost is chosen as
\begin{equation}
	\Phi = ||\theta(t_f)-\theta^*||^2/2+||\dot\theta(t_f)||^2/2,
	\nonumber
\end{equation}
so that the robotic arm stops at $t_f>0$ at the desired position~$\theta^*$. This linear problem can be solved with linear quadratic Gaussian (LQG) control \cite{Stengel}, however here a semi-analytical solution is obtained via implicit sampling.

\begin{figure}[tb]
\begin{center}
{\includegraphics[width=1\textwidth]{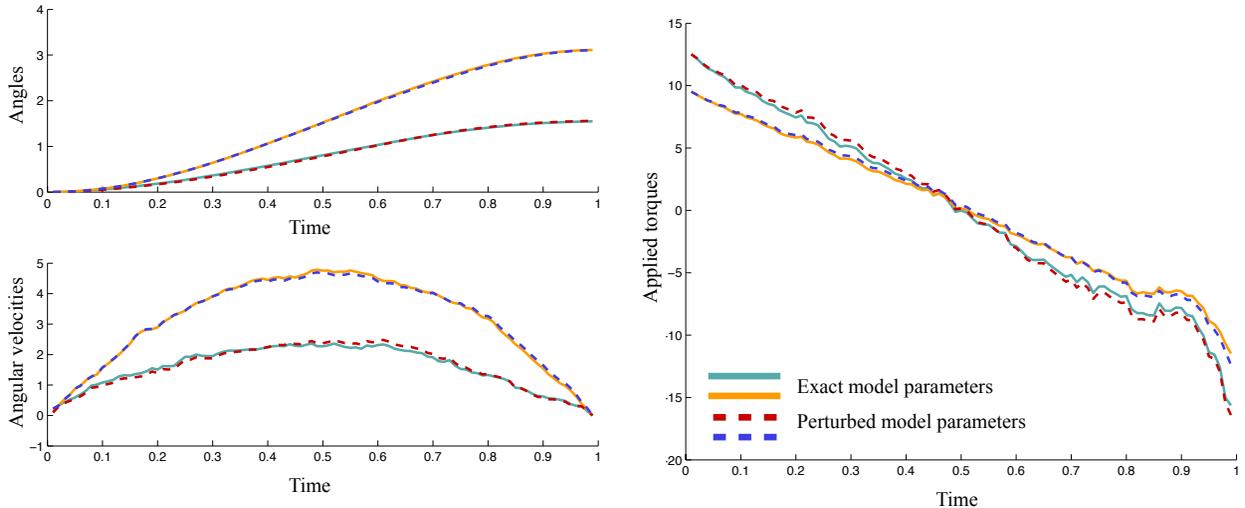}}
\caption{Simulation of a robotic arm under path integral control. 
Upper left: trajectories of angles $\theta_1$ (turquoise\slash red) and $\theta_2$ (orange\slash blue);
lower left: trajectories of angular velocities $\dot{\theta}_2$ (turquoise\slash red) and $\dot{\theta}_2$ (orange\slash blue);
right: applied torques $\tau_1$ (turquoise\slash red)and $\tau_2$ (orange\slash blue).
Solid turquoise and orange lines: trajectories and torques under optimal control with exact model parameters;
dashed red and blue lines: trajectories and torques under path integral control with false model parameters.}
\label{fig:RobotArmResults} 
\end{center}
\end{figure} 
As before, the time is discretized using a constant time step $\Delta t=(t_f-t)/n$. The discretized dynamics are 
\begin{eqnarray}
	z_{i+1}&=&z_i+ y_i\,\Delta t,\nonumber \\
	y_{i+1}&=&y_i+\Delta B,	\nonumber
\end{eqnarray}
where $y_i$, $z_i$ are 2-dimensional vectors whose elements are the discretized angular velocities ($\dot\theta$) and the discretized angles ($\theta$); note that only one of the above equations is driven by noise ($\Delta B$), since there is no reason to inject noise into the (trivial) equation $\dot\theta=\dot\theta$. While the noise propagates via the coupling to all variables, the pdfs that define the path in (\ref{eq:Integral}) are in terms of $y_n$s only.

The function $F$ in implicit sampling is thus a function of $y$ only,
\begin{equation}
	F=\frac{1}{2\gamma}y_n^2+\frac{1}{2\gamma}(z_n-\theta^*)^2+\frac{1}{2\Delta t}\sum_{i=1}^n(y_i-y_{i-1})^2,	\nonumber
\end{equation}
where $\gamma=r$ is chosen to satisfy (\ref{eq:NoiseCostRequirement}).
Because $F$ is quadratic the minimization is straightforward, and the algebraic equation (\ref{eq:ImplicitSamplingEquation}) can be solved with a linear map whose Jacobian is constant and given by the determinant of a Cholesky factor $L$ of the Hessian of~$F$ at the minimum (see above). Thus, equation~(\ref{eq:ImplicitSamplingSolution}) becomes
\begin{equation}
	\psi(x,t) = \frac{\exp(-\phi)}{(\Delta t)^{n}\det L}.\nonumber
\end{equation}
There is no need for generating samples, since the expected value is computed explicitly
(i.e.~by evaluating the integral analytically).

The robustness of the implicit sampling based path integral control is tested in numerical experiments.
To simulate that the ``real'' robotic arm behaves differently from the numerical model that
the path integral controller uses to find a control, one can give the controller false information about the 
parameters of the numerical model.
Here the parameters $m_1$ and $m_2$ are perturbed values of the true parameters
simulating that the robotic arm picked up a payload (of unknown weight).
Thus, the controller works with values $m_1=1.4$, $m_2=1$, whereas the ``true'' robotic arm
has parameters as in table \ref{tab:RoboticArm}.
The results of a simulation are shown in figure \ref{fig:RobotArmResults}
where state trajectories and controls are shown for a truly optimal controller 
(i.e.~one who has access to the exact model parameters),
and for the path integral controller which uses false model parameters.
One observes that the path integral controller can perform the desired task 
(moving the arm to a new position) and that its control and the resulting 
state trajectories closely follow those that are generated by a truly optimal controller.
This example thus indicates that the path integral controller can perform reliably and quickly 
in applications where model error may be an issue.

\subsection{Optimization via stochastic control}
One can set up a conventional optimization problem, i.e.~find $\min f(x)$ for a given smooth function $f$, as a stochastic optimal control problem as follows: find a control $u$ to minimize the cost function
\begin{equation}
\label{eq:StochOptimizaitonCost}
	C(x,t_f)=E\left[f(x(t_f))+\int_0^{t_f}u(x,t)^TRu(x,t)dt\right],
\end{equation}	
where the expectation is over trajectories of the SDE
\begin{equation}
	dx=u\,dt+\sqrt{\sigma}\,dW.\nonumber
\end{equation}
The idea is to make use of the stochastic component to explore valleys in the function $f$. The parameters one can tune are the noise level $\sigma$, the final time $t_f$ and the control cost $R$. 

Implicit sampling for this control problem requires minimization of the function
\begin{equation}
	F=\frac{1}{2\sigma\Delta t}\sum_{i=1}^n(y_i-y_{i-1})^2+\frac{1}{\gamma}f(y_n),\nonumber
\end{equation}
where the discretization is done using a constant time step $\Delta t=t_f/n$ as before. Note that the control approach to this problem inserts a quadratic term through which the space is (randomly) explored. 
\begin{figure}[tb]
\begin{center}
{\includegraphics[width=0.6\textwidth]{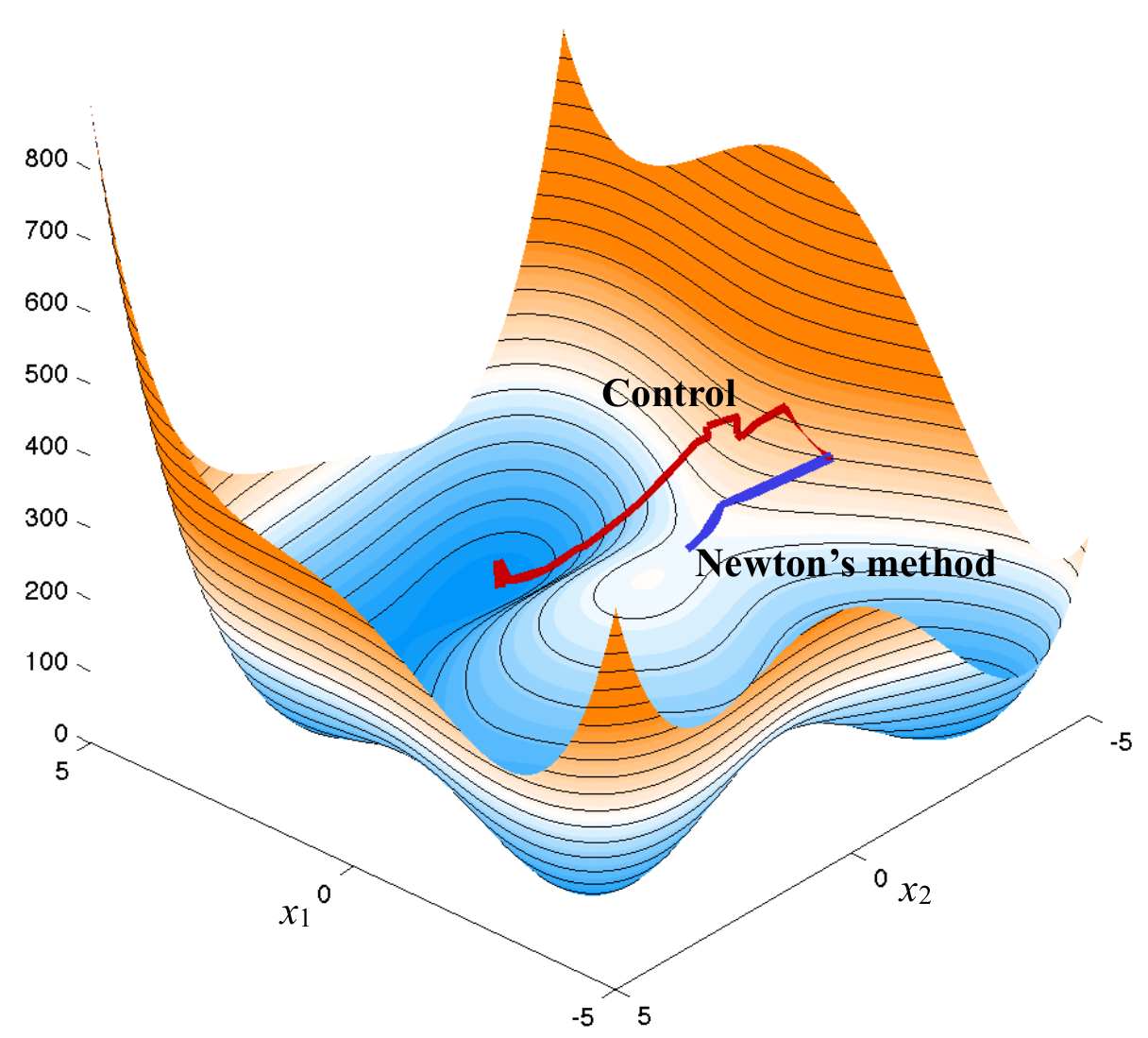}}
\caption{Newton's method and stochastic optimization for the Himmelblau function.} 
\label{fig:Himmelblau} 
\end{center}
\end{figure} 

These ideas are applied to minimize the Himmelblau function 
\begin{equation}
	f(x_1,x_2)=(x_1^2+x_2-11)^2+(x_1+x_2^2-7)^2,\nonumber
\end{equation}
which is a popular test of the performance of optimization algorithms because of its many local minima \cite{Himmelblau}. The parameters are $R=0.01,\sigma=0.01,\Delta t=1,t_f=20$, and the minimization is initialized at $(-1,-4)$. Figure \ref{fig:Himmelblau} shows the iterations obtained with implicit sampling and 50 samples, and, for comparison, the iterations of Newton's method.

In this example, the stochastic approach is successful and finds a much lower minimum than Newton's method. However, since each run of the stochastic approach is random, one may find another local minimum in another run. This could be useful when one needs to explore valleys or if one suspects the existence of other local minima. 

To test the reliability of the stochastic control approach, 70 experiments were performed, each starting at the same initial condition. The iterations of 5 such experiments are shown in table \ref{tab:StochasticOptimization}.
\begin{table}[tb]
\caption{Newton's method and 5 runs of controlled optimization}
\begin{center}
\begin{tabular}{ c c c c c c c }
Iteration & Newton & \multicolumn{5}{c}{Controlled Optimization}\\
\hline
0	&260.00	&260.00	&260.00	&260.00	&260.00	&260.00	\\
1	&195.36	&268.57	&266.47	&241.73 &337.77	&257.55	\\
2	&180.19	&258.45	&242.92	&239.19	&347.24	&263.20\\	
3	&178.43	&259.75	&268.95	&258.97 &360.16	&243.92\\
4	&178.34	&257.44	&279.41	&249.28	&373.97	&218.95\\
5	&178.34	&214.25 &243.23	&224.21	&336.47	&201.69\\ 
6	&178.34	&184.70	&320.29	&256.08	&476.23	&179.34\\
7	&178.34	&161.31	&331.32	&230.16	&394.82	&157.00 \\
8	&-		&147.39	&330.16	&218.89	&316.97	&150.73\\
9	&-		&144.04	&255.65	&175.61	&245.24	&150.33\\
10	&-		&101.05	&221.05	&180.92	&176.98	&144.97\\
11	&-		&76.729	&188.49	&174.13	&180.53	&105.48\\
12	&-		&48.46	&151.54	&188.66	&131.86	&91.42 \\
13	&-		&33.05  &115.34	&171.75	&101.82	&67.00\\
14	&-		&29.90	&87.37	&147.43	&56.51	&41.29\\
15  &-		&41.34	&62.89	&104.69	&35.11	&10.65\\
16	&-		&31.03	&43.69	&72.71	&31.48	&6.08\\
17	&-		&18.12	&22.62	&42.44	&14.53 	&1.79\\
18	&-		&1.83	&9.195	&30.55	&9.52	&2.61	\\
19	&-		&0.46	&1.90	&0.54	&0.97	&1.24	\\
\hline
\end{tabular}
\end{center}
\label{tab:StochasticOptimization}
\end{table}
All 70 runs ended up close to the local minimum around $x_1=2.5$, $x_2=-2.5$. The average value of $F$ is an approximation to $C$ in (\ref{eq:StochOptimizaitonCost}) and, with 70 experiments was found to be $C\approx 1.75$, which corresponds to a much lower value of $f$ than the value $f_{\min}^{\text{Newton}} = 178.34$, that one obtains with Newton's method.

\section{Application to Monte Carlo localization and SLAM}\label{sec:MCLandSLAM}
Consider a mobile robot that navigates autonomously and, as it moves, collects noisy measurements about its motion as well as scans of its environment. If the scans reveal locations of features that are known to the robot, i.e.~if the robot has a map of its environment, then it can localize itself on this map. This problem is known as localization. When the robot has no map of its environment, then it must construct a map while localizing itself on it, leading to the problem of simultaneous localization and mapping (SLAM)  \cite{CaseStudyBook,SLAM1,RobotBook}. Efficient solutions to the localization and SLAM problems are a fundamental requirement for autonomous robots, which must move through unknown environments where no global positioning data are available, for example indoors, in abandoned mines, underwater, or on far-away planets \cite{SCIENCE,ThrunMines,ThrunBook}. Application of implicit sampling to localization and SLAM is the subject of the next two sections, where the algorithms will also be tested using experimental data obtained by Nebot and colleagues at the University of Sydney \cite{UTE_DataSet}.

\subsection{Monte Carlo localization}
The localization problem can be formulated as follows. Information about the initial state of the robot is available in the form of a pdf $p(x_0)$, where $x_0$ is an $m$-dimensional vector whose elements are the state variables (e.g.~position and velocity of the robot). A probabilistic motion model defines the pdf
\begin{equation}
\label{eq:MotionModel} 
p(x_n|x_{n-1},u_{n}),
\end{equation}
where $n=1,2,\dots$ is discrete time and where $u_n$ is a $p$-dimensional vector of known ``controls'', e.g.~the odometry. A data equation describes the robot's sensors and defines the pdf
\begin{equation}
\label{eq:ObsModel}
	p(z_n|x_n,u_n,\Theta),
\end{equation}
where $z_n$ is a $k$-dimensional vector whose elements are the data ($k\leq m$) and where $\Theta$ is the map. For example, one can use the ``landmark model'', for which the map describes the coordinates of a collection of obstacles, called``landmarks''; the data are measurements of range and bearing of the position of the robot relative to a landmark. The landmark model and range and bearing data will be used in the applications below, but the algorithms described are more generally applicable. 

The motion and sensor model jointly define the conditional pdf $p(x_{0:n}|z_{0:n},u_{0:n},\Theta)$, where the short-hand notation $x_{0:n }=\{x_0,x_1,\dots,x_n\}$ is used to denote a sequence of vectors. By using Bayes' rule repeatedly, one can derive the recursive formula 
\begin{equation*}
	p(x_{0:n}|z_{0:n},u_{0:n},\Theta) =p(x_{0:n-1}|z_{0:n-1},u_{0:n-1},\Theta)\frac{p(x_n|x_{n-1},u_n)p(z_n|x_n,u_n,\Theta)}{p(z_n|z_{0:n-1})}.
\end{equation*}
Approximations of this pdf are used to localize the robot. For example, methods based on the Kalman filter \cite{Kalman1960} (KF) construct a Gaussian approximation which is often problematic because of nonlinearities in the model and data. Moreover, KF-based localization algorithms can diverge, for example, in multi-modal scenarios (i.e.~if the data are ambiguous) or during re-localization after system failure \cite{FastSLAM}. The basic idea of Monte Carlo localization (MCL) is to relax the Gaussian assumptions required by KF, and to apply importance sampling to the localization problem \cite{ThrunMCL}. The method proceeds as follows. 

Given $M$ weighted samples $\{X_{0:n-1}^j,W_{n-1}^j\}$, $j=1,\dots,M$, which form an empirical estimate of the pdf $p(x_{0:n-1}|z_{0:n-1},u_{0:n-1},\Theta)$
 at time $n-1$, one updates each particle to time $t_n$ given the datum~$z_n$. In standard MCL \cite{ThrunMCL}, this is done by choosing the importance function 
\begin{equation*}
	(p_0)_n \propto p(x_{0:n-1}|z_{0:n-1},u_{0:n-1},\Theta) p(x_n|X^j_{n-1},u_n),
\end{equation*}
i.e.~the robot location is predicted with the motion model and then this prediction, $X^j_n$, is weighted by 
\begin{equation*}
	W^j_n \propto W^j_{n-1} p(z_n|X^j_n,u_n,\Theta),
\end{equation*}
to assess how well it fits the data. However, this choice of importance function can be inefficient, especially if the motion model is noisy and the data are accurate (as is often the case \cite{ThrunBook}). The reason is that many of the samples generated by the model are unlikely with respect to the data and the computations used to generate these samples are wasted.

Implicit sampling can be used to speed up the computations. This requires in particular that implicit sampling, as described in section 2, is applied to the $M$ functions
\begin{equation}
\label{eq:F_IS_MCL}
	F^j(x_n) = -\log\left(p(x_n|X^j_{n-1},u_n)p(z_n|x_n,u_n,\Theta)\right).
\end{equation}
One needs $M$ functions $F^j$, $j=1,\dots,M$, one per sample, because the recursive problem formulation requires a factorization of the importance function (compare to sections 2 and 3, where only one function was needed). Here each function $F_j$ is parametrized by the location of the sample at time $n-1$, $X^j_{n-1}$, the control $u_n$ and the datum $z_n$; the variables of these functions are the location of the robot at $t_n$, $x_n$.

For MCL with implicit sampling, each function $F^j$ must be minimized, e.g.~using Newton's method. After the minimization, one can generate samples by solving the stochastic equation (\ref{eq:ImplicitSamplingEquation}) with the techniques described in section 2. 
Using information from derivatives in sampling for MCL has also been considered 
in \cite{Biswas2011,Biswas2013}.

\subsubsection{Implementation for the University Car Park data set} \label{sec:MCL_Example}
The University Car Park data set \cite{UTE_DataSet} is used to demonstrate the efficiency of MCL with implicit sampling. The scenario is as follows. A robot is moving around a parking lot and steering and speed data are collected via a wheel encoder on the rear left wheel and a velocity encoder. An outdoor laser (SICK LMS 221) is mounted on the front bull-bar and is directed forward and horizontal (see \cite{UTE_Manual} for more information on the hardware). The laser returns measurements of relative range and bearing to different features. Speed and steering data are the controls of a kinematic model of the robot and the laser data are used to update this model's predictions about the state of the robot. 

The motion model is the forward Euler discretization of the continuous time 2D-model in \cite{UTE_Manual},
\begin{equation*}
	x_{n+1}=x_n+R(x_n,u_n)\delta+\Delta B_n, \quad \Delta B_n=\mathcal{N}(0,(\Sigma_1)_n),
\end{equation*}
where $R(x_m,u_n)$ is a $3$-dimensional vector function, $\delta$ is the time step, and $(\Sigma_1)_n$ is a given covariance matrix.  
The data equation is
\begin{equation*}
z_{n} = h(x_n,\Theta)+V_n,\quad V_n \sim \mathcal{N}(0,\Sigma_2),
\end{equation*}
where $h$ is a 2-dimensional vector function that maps the position of the robot to relative range and bearing and where $\Sigma_2$ is a $2 \times 2$ diagonal matrix (see \cite{UTE_Manual} and the appendix for the various model parameters).

With the above equations for motion model and data, the functions $F^j$ of implicit sampling in (\ref{eq:F_IS_MCL})  become
\begin{equation}
\label{eq:MCL_F}
F^j(x_n) = (x_n-f^j_n)^T\left((\Sigma_1)_{n}\right)^{-1}(x_n-f^j_n)+\displaystyle \sum_{i=1}^p (h(x_n)-z_n^i)^T\Sigma_2^{-1}(h(x_n)-z_n^i),
\end{equation}
where $f^j_n =X^j_{n-1}+R(X^j_{n-1},u_n)\delta$ and where $z^i_n$ denotes the measurement for the $i$th landmark at time $n$. These $F^j$s are minimized with Newton's method and samples are generated by solving a quadratic equation as explained in section~2. The gradient and Hessian in Newton's method were computed analytically (see appendix). If uneven weights were observed, a ``resampling'' was done using Algorithm~2 in~\cite{GordonReview}. Resampling replaces samples with a small weight with samples with a larger weight without introducing significant bias. If no laser data are available, all samples evolve according to the model equation (\ref{eq:MotionModel}).

\begin{figure}[tb]
\begin{center}
{\includegraphics[width=1\textwidth]{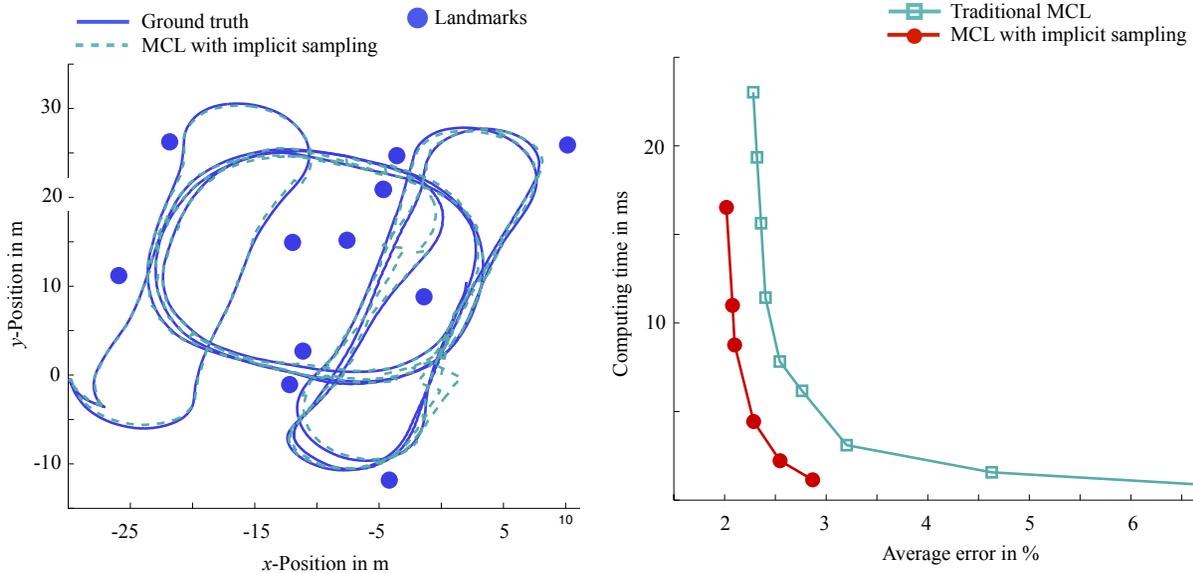}}
\caption{Simulation results of MCL with implicit sampling and standard MCL. Left: the path of the robot and its reconstruction via MCL with implicit sampling (100 samples). Right: average computing time as a function of the average error for standard MCL and MCL with implicit sampling.} 
\label{fig:MCL_Results} 
\end{center}
\end{figure}

The results of MCL with implicit sampling are shown in the left panel of Figure \ref{fig:MCL_Results}. 
The ``ground truth'' shown here is the result of an MCL run with GPS data (these GPS data however are not used in implicit sampling); the large blue dots are the positions of the landmarks. The reconstruction by the MCL algorithm (dashed turquoise line) is an approximation of the conditional mean, which is obtained by averaging the samples. One observes that MCL with implicit sampling gives accurate estimates whenever laser data are available. After periods during which no laser data were available, one can observe a strong ``pull'' towards the true state trajectory, as soon as data become available, as is indicated by the jumps in the trajectory e.g.~around $x=0$, $y=0$. 

The efficiency of MCL with implicit sampling is studied by running the algorithm with $10,\, 20,\, 40,\, 80,\, 100$ and $150$ samples, i.e.~with increasing computational cost, and comparing the reconstruction of the MCL with the true path. The error is measured by the Euclidean norm of the difference between the ``ground truth'' (see above) and the reconstruction by the MCL algorithm (scaled by the Euclidean norm of the ground truth). The standard MCL method was also applied with $10,\, 20,\, 40,\, 80,\, 100,\,150,\,250$ and $300$ samples. The results are shown in table \ref{tab:MCLResults}.
\begin{table}[tb]
\caption{Simulation results for MCL}
\begin{center}
\begin{tabular}{c c c c c}
& \multicolumn{2}{c}{Implicit sampling} & \multicolumn{2}{c}{Standard sampling}\\
\hline
\# of samples  & Error & CPU time & Error & CUP time \\
  & in \% & in ms & in \% & in ms \\
  10	& 2.87&1.14 & 6.91 &0.81\\
20	& 2.55&2.23 & 4.63 &1.57 \\
40	&2.29&4.42 & 3.20 & 3.12 \\
80	& 2.10&8.80 & 2.76 & 6.20\\
100	&2.08  &11.0 & 2.54 &7.75\\
150	&2.02  &16.5 & 2.40 &11.6\\
200	& - & -   & 2.36 &15.8\\
250	&- & -   & 2.32 & 19.6 \\
300	& -& -   & 2.28 &  23.4 \\
\hline
\end{tabular}
\end{center}
\label{tab:MCLResults}
\end{table}

It is clear that the both MCL algorithms converge to the same error, since both algorithms converge to the conditional mean as the number of samples (and hence, the computational cost) goes to infinity. However the convergence rate of MCL with implicit sampling is faster, which can be seen from the right panel of figure \ref{fig:MCL_Results}, where the computing time is plotted as a function of the error. In this figure, the area between the lines corresponds to the improvement of implicit sampling over the standard method, and it is clear that implicit sampling is more efficient than the standard method. The improvement is particularly pronounced if a high accuracy is required, in which case MCL with implicit sampling can be orders of magnitudes faster than the standard method.

\subsection{Implicit sampling for online SLAM} \label{sec:SLAM}
In SLAM, one is given the probabilistic motion model~(\ref{eq:MotionModel}) and a data equation~(\ref{eq:ObsModel}), and one needs to estimate the position of the robot as well as the configuration and geometry of the map; the conditional pdf of interest is thus
\begin{equation}
\label{eq:SLAMpdf}
	p(x_{0:n},\Theta|z_{0:n},u_{0:n},\eta_{0:n}),
\end{equation}
where the map, $\Theta$, is a variable (not a given parameter as in MCL); the variable $\eta_{0:n}$ is the ``data association'', i.e.~it maps the data to the known features, or creates a new feature if the data are incompatible with current features \cite{SmithCheesman,Durrant-Whyte}. Here it is assumed that the data association is done by another algorithm, so that $\eta_{0:n}$ in~(\ref{eq:SLAMpdf}) can be assumed to be given (in the numerical implementation in section \ref{sec:SLAM_Example}, a maximum-likelihood algorithm is used \cite{FastSLAM}).

Note that the SLAM state-vector is different from the state of the localization problem: it is the number of variables needed to describe the robot's position, $x$, and the coordinates of the features of the map, $\Theta$. If the map contains many features (which is typically the case), then the state dimension is large and this makes KF based SLAM prohibitively expensive. The reason is that KF SLAM requires dense matrix operations on matrices of size $N$ (where $N$ is the number of features), due to correlations between the robot's position and the features \cite{SmithCheesman}. KF SLAM thus has $N^2$ memory requirements which make it infeasible for realistic $N$.

Monte Carlo approaches to SLAM reduce the computational cost by exploring conditional independencies \cite{Thrun1998}. In particular, the various features of the map conditioned on the robot path are independent, which implies the factorization 
\begin{equation*}
	p(\Theta|x_{0:n},z_{0:n},u_{0:n},\eta_{0:n}) = \displaystyle\prod_{k=1}^N p(\theta^k|x_{0:n},z_{0:n},u_{0:n},\eta_{0:n}),
\end{equation*}
where $\theta^k$, $k=1,\dots,N$ are the features of the map~\cite{Murphy,FastSLAM}. This factorization makes it possible to update one feature at a time and thus reduces memory requirements. Moreover, ``online'' SLAM will be considered here, i.e.~the map and the robot's position are constructed recursively as data become available, using
\begin{eqnarray}
	p(x_{0:n},\Theta|z_{0:n},u_{0:n},\eta_{0:n}) &\propto& p(x_{0:n-1},\Theta|z_{0:n-1},u_{0:n-1},\eta_{0:n-1}) \nonumber \\
	&& \times \, p(x_n|x_{n-1}, u_n)p(z_n|\theta_n,x_n,\eta_n), \label{eq:RecursiveSLAMpdf}
\end{eqnarray}
where, $\theta_n$ is the feature observed at time $n$. Alternatively, one can wait and collect all the data, and then assimilate this large data set in one sweep; such a ``smoothing'' approach is related to graph based SLAM (see e.g.~\cite{GraphSLAM}) but will not be discussed further in this paper.

For online probabilistic SLAM, one assumes that $M$ weighted samples $\{X^j_{0:n-1},\Theta^j,W^j_{n-1}\}$ approximate 
\begin{equation*}
	p(x_{0:n-1},\Theta|z_{0:n-1},u_{0:n-1},\eta_{0:n-1}),
\end{equation*}	
at time $n-1$ and then updates the samples to time $n$ by applying importance sampling to
\begin{equation*}
p(x_n|x_{n-1}, u_n)p(z_n|\theta_n,x_n,\eta_n).
\end{equation*}
Here it is assumed that only one feature is observed at a time (which is realistic \cite{FastSLAM}), however the extension to observing multiple features simultaneously is straightforward. 

Various importance functions have been considered and tested in the literature. For example, the fastSLAM algorithm \cite{Montemerlo2002} was shown to be problematic in many cases due to the poor distribution of its samples \cite{Montemerlo03a,FastSLAM}. The reason is that the samples are generated by the motion model and, thus, are not informed by the data. As a result, the samples are often unlikely with respect to the data. The fastSLAM~2.0 algorithm \cite{Montemerlo03a,FastSLAM} ameliorates this problem by making use of an importance function that depends on the data. The algorithm was shown to outperform fastSLAM and was proven to converge (as the number of samples goes to infinity) to the true SLAM posterior for linear models. The convergence for nonlinear models is not well understood (because of the use of extended Kalman Filters (EKF) to track and construct the map). Here a new SLAM algorithm with implicit sampling is described. The algorithm converges for nonlinear models at a computational cost that is comparable to the cost of fastSLAM 2.0. 

Online SLAM with implicit sampling requires that implicit sampling is applied to the $M$ functions $F^j$ (one per sample), whose variables are the position at time $t_n$, $x_n$, and the location of the feature $\theta_n$, observed at time $t_n$:
\begin{equation*}
	F^j(x_n,\theta)= -\log p(x_n|X^j_{n-1}, u_n)p(z_n|\theta^j_n,x_n,\eta_n^j);
\end{equation*}
the position of the $j$th sample at time $t_{n-1}$ is a parameter. As in MCL, one needs $M$ functions here due to the recursive problem formulation. The minimization problems of implicit sampling are
\begin{equation*}
	\phi^j = \min_{x_n,\theta_n}F^j(x_n,\theta_n),
\end{equation*}
and samples $\{X^j_n,\theta^j_n\}$ are generated by solving the stochastic equations
(\ref{eq:ImplicitSamplingEquation}). One can use the same technique for solving these equations as described in section 2. Moreover, if the feature $\theta_n$ is already known, the SLAM algorithm reduces to the MCL algorithm with implicit sampling described above. 

Note that the memory requirements of SLAM with implicit sampling are linear in the number of features, because the algorithm makes use of the conditional independence of the features given the robot path, so that, at each step, only one feature needs to be considered. Moreover, SLAM with implicit sampling converges to the true SLAM posterior (under the usual assumptions about the supports of the importance function and target pdf \cite{Geweke}) as the number of samples goes to infinity; this convergence is independent of linearity assumptions about the model and data equations. Convergence for fastSLAM and fastSLAM~2.0 however can only be proven for linear Gaussian models \cite{FastSLAM}. 

\subsubsection{Numerical experiments} \label{sec:SLAM_Example}
The applicability and efficiency of SLAM with implicit sampling is demonstrated with numerical experiments that mimic the University Car Park data set \cite{UTE_DataSet}. Here speed and steering data of \cite{UTE_DataSet} are used, however the laser data are replaced by synthetic data, because the data in \cite{UTE_DataSet} are too sparse (in time) for successful online SLAM (even EKF SLAM, often viewed as a benchmark, does not give satisfactory results).

\begin{figure}[tb]
\begin{center}
{\includegraphics[width=1\textwidth]{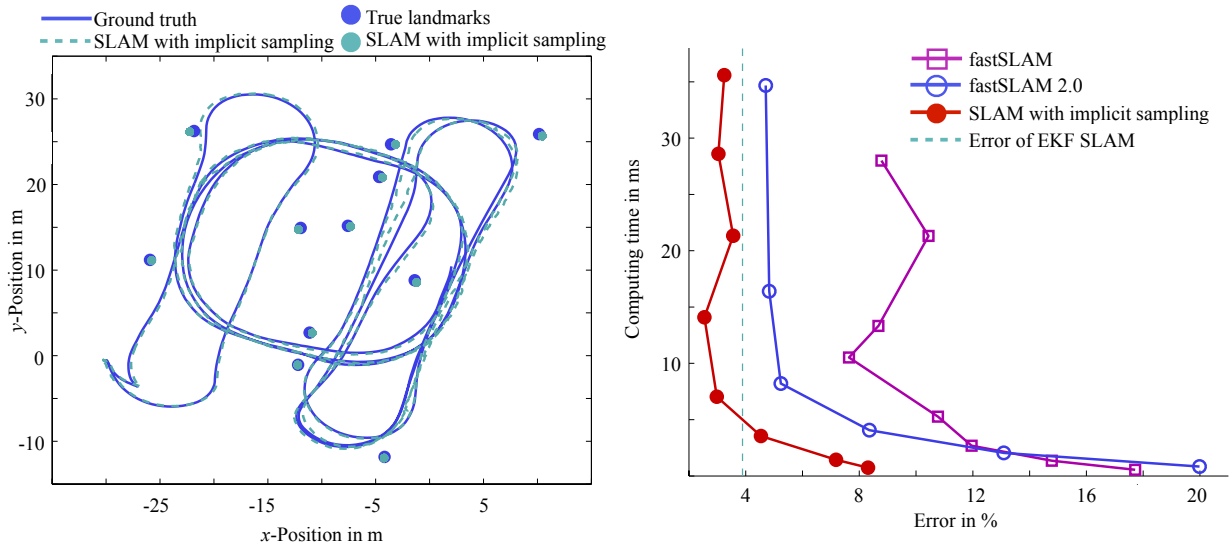}}
\caption{Simulation results for SLAM. Left: ground truth and reconstruction via SLAM with implicit sampling (100 samples). Right: computing time as a function of the average error of fastSLAM, fastSLAM~2.0 and SLAM with implicit sampling.} 
\label{fig:SLAM_Results} 
\end{center}
\end{figure}

In these synthetic data experiments, a virtual laser scanner returns noisy range and bearing measurements of features in a half-circle with $15$m radius. If the laser detects more than one feature, one is selected at random and returned to the SLAM algorithm. Since each synthetic data set is a random event, 50 synthetic data sets are used to compute the average errors for various SLAM algorithms. These errors are defined as the norm of the difference of the ground truth and the conditional mean computed with a SLAM algorithm. The performances of EKF SLAM, fastSLAM, fastSLAM~2.0 and SLAM with implicit sampling are compared using these synthetic data sets. All SLAM methods make use of the same maximum likelihood algorithm for the data association. A Newton method was implemented for the optimization in implicit sampling, and the quadratic approximation of $F$ was used to solve the algebraic equations (\ref{eq:ImplicitSamplingEquation}). A typical outcome of a numerical experiment is shown in the left panel of figure \ref{fig:SLAM_Results}, where the true path and true positions of the landmarks as well as their reconstructions via SLAM with implicit sampling (100 samples) are plotted. 

The results of 50 numerical experiments with fastSLAM, fastSLAM 2.0 and SLAM with implicit sampling are shown in table \ref{fig:SLAM_Results}. 
\begin{table}[tb]
\caption{Simulation results for SLAM}
\begin{center}
\begin{tabular}{c c c c c c c}
 & \multicolumn{2}{c}{fastSLAM} & \multicolumn{2}{c}{fastSLAM 2.0} & \multicolumn{2}{c}{Implicit sampling}\\
 \hline
Samples & CPU time & Error & CPU time & Error & CPU time & Error\\
  & in ms&  in \% &  in ms & in \%&  in ms &  in \%\\
10  &  -  & -   &  -  & -   & 7.31& 8.31\\
20	& 5.55& 17.7& 8.36& 20.0& 14.4& 7.18 \\
50	& 13.6& 14.8& 20.5& 13.1& 35.4& 4.54\\
100	& 26.9& 12.0& 40.7& 8.36& 70.3& 2.98\\ 
200	& 52.7& 10.8& 82.1& 5.24& 140.9& 2.55\\ 
400	& 105.2& 7.64& 164.0& 4.83& 286.0& 3.25\\ 
500	& 133.2& 8.67& 346.7& 4.71& 355.9& 3.25\\ 
800	& 213.2& 10.4& -    &-    &-     &-\\ 
1000& 280.1& 8.78& -    &-    &-     &-\\ 
\hline
\end{tabular}
\end{center}
\label{tab:SLAMResults}
\end{table}
It was observed in this example that EKF SLAM gives an average error of $3.89\%$ at an average computation time of $0.43$ ms and, thus, is computationally efficient and accurate. However, computational efficiency disappears in scenarios with larger maps due to the $O(N^2)$ memory requirements ($N$ is the number of features in the map). The error of fastSLAM~2.0 (whose memory requirements scale linearly with $N$) decreases with the number of samples, and it seems as if the fastSLAM~2.0 solution converges (as the number of samples and, thus, the computation time increases) to the EKF solution. This is intuitive because fastSLAM~2.0 uses EKFs to track the map. The fastSLAM algorithm, on the other hand, converges to an error that is larger than in EKF SLAM or fastSLAM~2.0. SLAM with implicit sampling converges to an error that is smaller than in EKF SLAM. The reason is as follows: implicit sampling requires no linearity assumptions, and, therefore, the empirical estimate it produces converges (as the number of samples goes to infinity) to the true SLAM posterior. 

Moreover, the convergence rate of SLAM with implicit sampling is faster than for any other SLAM method, due to the data informed importance function. Thus, the approximation of the conditional mean (the minimum mean square error estimator) one obtains with implicit sampling is accurate even if the number of samples is relatively small. As a result, SLAM with implicit sampling is more efficient than the other SLAM algorithms considered here. This is illustrated in the right panel of figure \ref{fig:SLAM_Results}, where the computing time is shown as a function of the average error. As in MCL (see above), the area between the various curves indicates the improvement in efficiency. Here implicit sampling is the most efficient method, giving a high accuracy (small error) at a small computational cost. Moreover, SLAM with implicit sampling can achieve an accuracy which is higher than the accuracy of all other methods (including EKF SLAM),
because it is a convergent algorithm.

\section{Conclusions}
Implicit sampling is a Monte Carlo scheme that localizes the high-probability region of the sample space via numerical optimization, and then produces samples in this region by solving algebraic equations with a stochastic right hand side. The computational cost per sample in implicit sampling is larger than in many other Monte Carlo sampling schemes that randomly explore the sample space. However, the minimization directs the computational power towards the relevant region of the sample space so that implicit sampling often requires fewer samples than other MC methods. There is a trade-off between the cost-per-sample and the number of samples and this trade-off was studied in this paper in the context of three applications: Monte Carlo localization (MCL) and probabilistic online SLAM in robotics, as well as path integral control. 

In path integral control one solves the stochastic Hamilton-Jacobi-Bellman (HJB) equation with a Monte Carlo solver. This has the advantage that no grid on the state space is needed and the method is therefore (in principle) applicable to control problems of relatively large dimension. Implicit sampling was applied in this context to speed up the Monte Carlo calculations. The applicability of this approach was demonstrated on two test problems. Path integral control was also used to find local minima in non-convex optimization problems. 

An implementation of implicit sampling for MCL was tested and it was found that implicit sampling performs better than standard MCL (in terms of computing time and accuracy), especially if the data are accurate and the motion model is noisy (which is the case typically encountered in practice \cite{ThrunBook}). Similarly, implicit sampling for SLAM outperformed standard algorithms (fastSLAM and fastSLAM~2.0). Under mild assumptions, but for nonlinear models, SLAM with implicit sampling converges to the true SLAM posterior as the number of samples goes to infinity, at a computational cost that is linear in the number of features of the map.

\section*{Acknowledgements}
I thank Professor Alexandre J.~Chorin of UC Berkeley for many interesting technical discussions and for bringing path integral control to my attention. I thank Dr.~Robert Saye of Lawrence Berkeley National Laboratory for help with proofreading this manuscript. 
This material is based upon work supported by the 
U.S.~Department of Energy, Office of Science,
Office of Advanced Scientific Computing Research, Applied Mathematics program under contract DE-AC02005CH11231, 
and by the National Science Foundation under grant DMS-1217065.

\newpage
\section*{Appendix}
\begin{figure}[tb]
\begin{center}
{\includegraphics[width=0.5\textwidth]{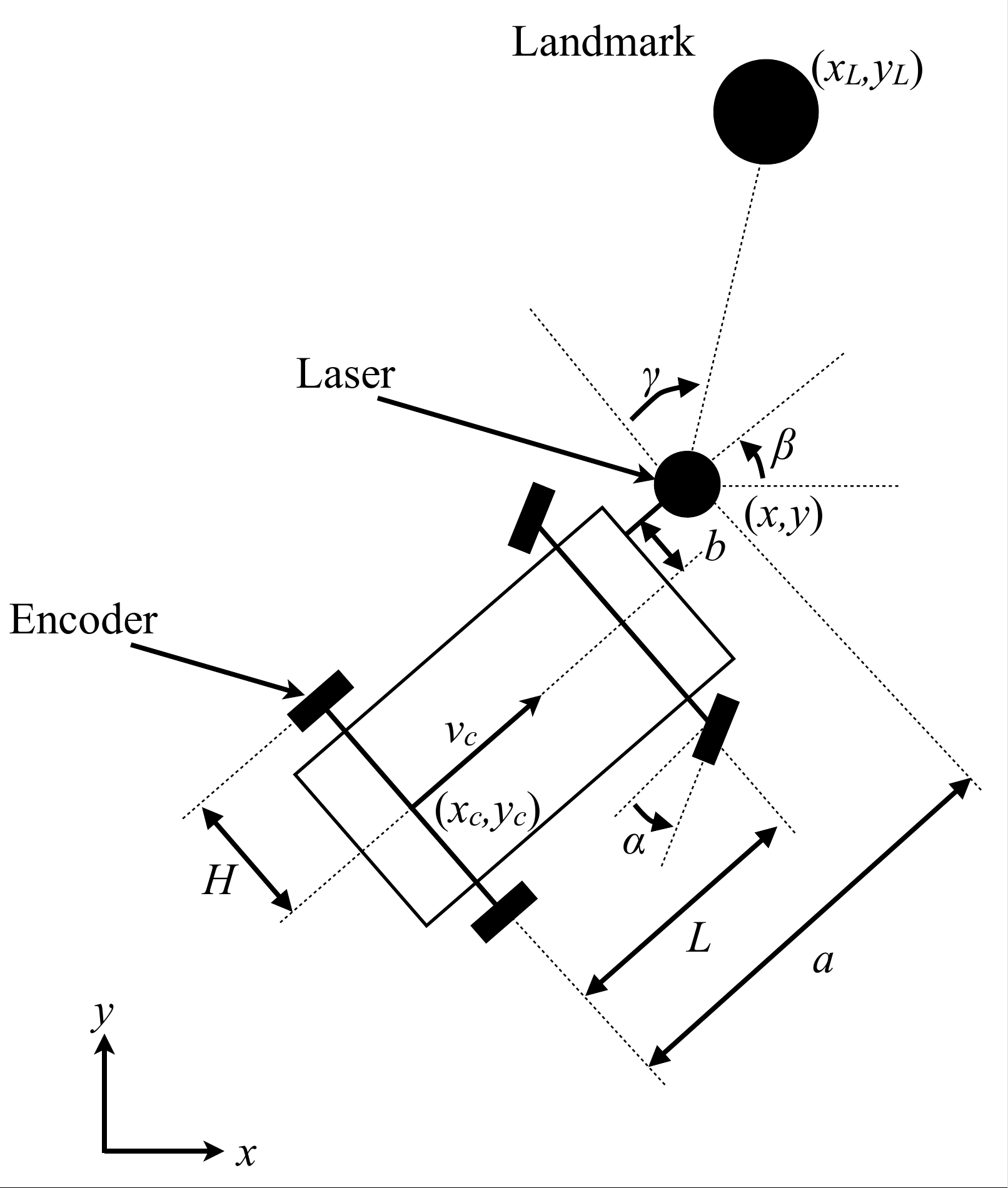}}
\caption{Vehicle kinematics (a version of Figure 1 in \cite{UTE_Slam}).} 
\label{fig:Kinematics} 
\end{center}
\end{figure}

A kinematic model of the robot is described in \cite{UTE_Slam} and shown for convenience in figure \ref{fig:Kinematics}. 
The robot is controlled by the speed $v_c$ and the steering angle~$\alpha$. It can be shown that the derivatives of the position and orientation of the axle of the robot are
\begin{equation*}
\frac{dx_c}{dt} = v_c \cos(\beta),\quad
\frac{dy_c}{dt} = v_c \sin(\beta),\quad
\frac{d\beta}{dt} = \frac{v_c}{L}\tan(\alpha) .
\end{equation*}
The velocity $v_l$ however is measured at the rear left wheel and, thus, must be translated to the axle:
\begin{equation*}
v_c=\frac{Lv_l}{L-\tan{\alpha}H}.
\end{equation*} 
The position of the laser can be obtained from the position of the axle by using
\begin{equation*}
x=x_c+a\cos(\beta)-b\sin(\beta),\quad y=y_c+a\sin(\beta)+b\cos(\beta),
\end{equation*}
so that the forward model (without noise) is 
\begin{eqnarray*}
\frac{dx}{dt}&=&v_c\left(\cos(\beta)-\tan{\alpha}\left(a\sin(\beta)-b\cos(\beta)\right)\right),\\
\frac{dy}{dt}&=&v_c\left(\sin(\beta)+\tan{\alpha}\left(a\cos(\beta)-b\sin(\beta)\right)\right),\\
\frac{d\beta}{dt}&=&\frac{v_c}{L}\tan(\alpha),
\end{eqnarray*}
which, in a more compact notation, can be written as
\begin{equation*} 
	\frac{dx}{dt}=R(x,u),
\end{equation*} 
where $x=(x,y,\beta)^T$ and $u=(v_c,\alpha)^T$. The uncertainty in the model is due to errors in the inputs and in the model \cite{UTE_Slam}: 
\begin{eqnarray*}
	dx&=&R(x,\hat{u}(t))\,dt+Q\,dW_x,\\
	   \hat{u}(t)\,dt&=&u(t)\,dt+ P\,dW_u.
\end{eqnarray*}
where $W_u$ and $W_x$ are 2-, respectively 3-dimensional vectors whose components are independent Brownian motions, $P$ and $Q$ are real $2\times 2$, respectively $3\times 3$, matrices and $u(t)$ is the unperturbed control signal. The above equations are linearized in $\hat{u}$ to obtain
\begin{equation*}
dx = R(x,u)\,dt+\frac{dR}{du}P\,dW_u+Q\,dW,
\end{equation*}
and the stochastic forward Euler method \cite{Kloeden} with time step $\delta$ is used to discretize the equations:
\begin{equation*}
	x^{n+1}=x^n+R(x^n,u^n)\delta+\Delta W^n, \quad \Delta W^n=\mathcal{N}(0,\Sigma_1^n),
\end{equation*}
where
\begin{equation*}
	\Sigma_1^n = \delta \left(\frac{dR}{du}PP^T\frac{dR}{du}^T+QQ^T\right),
\end{equation*}
and $x^n=x(n\delta),u^n=u(n\delta)$. Table \ref{tab:ModelParameters} lists the numerical values of all parameters of the probabilistic motion model.
\begin{table}[tb]
\caption{Model parameters.}
\begin{center}
\begin{tabular}{l c c}
Description &	Parameter& Value \\ \hline
Wheel base &	    $L$ & 2.83m \\
Width of the car &		$H$ & 0.76m\\
Horizontal offset of laser &		$b$ & 0.5m\\
Rear axis to laser &		$a$ & 3.78m\\ 
Time step	& 	$\delta$ & 0.025	\\
Covariance of model noise &	 	$Q$ &  $Q_{11}=0.1$, $Q_{22}=0.1$, $Q_{33}=0.1\pi/180$\\
& & $Q_{ij}=0$ if $i\neq j$\\
Covariance of noise in controls &		$P$ &  $P_{11}=P_{22}=0.5$, $P_{ij}=0$ if $i\neq j$\\
	\hline
\end{tabular}
\end{center}
\label{tab:ModelParameters}
\end{table}
The steering data $\alpha$ and velocity data $v_c$ are taken from the data set \cite{UTE_DataSet}.

The data are measurements of range and bearing of the features in the map $\theta$ relative to the robot and are obtained by the laser. Only ``high intensity points'' from the data set \cite{UTE_DataSet} are used. The laser is modeled by the data equation
\begin{equation*}
z^{j,n} = h(\hat{x}^n)+V^n,\quad V^n \sim \mathcal{N}(0,\Sigma_2),
\end{equation*}
where $\Sigma_2$ is a $2 \times 2$ diagonal matrix whose elements are $\sigma_1=0.05$ and $\sigma_2=0.05\pi/180$ and
$$
h(x^n)=\left(\begin{array}{c}
||x^n-m^j||\\
\mbox{atan2}\left(\frac{m_2^n-\hat{x}_2^n}{m_1-\hat{x}_1^n}\right)-x_3^n+\frac{\pi}{2},
\end{array}
\right),
$$
where $\hat{x}^n_j$ is the $j$th element of the vector $\hat{x}^n$, the ``true'' robot position, and $m^j_1$, $m^j_2$ are the $x$- and $y$-coordinates of the $j$th element of a landmark; the norm $||\cdot||$ is the Euclidean norm.

The gradient and Hessian of $F_j$ in (\ref{eq:MCL_F}) are:
\begin{align}
\nabla F_j &= \left(\Sigma_1^{n}\right)^{-1}(x-f^n) + \displaystyle \sum_{i=1}^p 
\left(\begin{array}{c}
\frac{x_1-m^i_1}{\sigma_1 r} (r_i-z^i_1)+\frac{x_2-m^i_2}{\sigma_2 r_i^2}(z^i_2-\hat\theta^i)\\
\frac{x_2-m^i_2}{\sigma_1 r}(r_i-z^i_1) -\frac{x_1-m^i_1}{\sigma_2 r_i^2}(z^i_2-\hat\theta^i)\\
\frac{z^i_2-\hat\theta^i}{\sigma_2}
\end{array}\right),\nonumber \\
H_j &=\left(\Sigma_1^{n}\right)^{-1}+ \displaystyle \sum_{i=1}^p\hat{H}_i,
\nonumber
\end{align}
where the elements of the $3\times 3$ matrices $\hat{H}^i$ are
\begin{align}
\hat{H}^i_{11} =& \frac{r_i-z^i_1}{\sigma_1 r_i}-\frac{(x_1-m^i_1)^2(r_i-z^i_1)}{\sigma_1 r_i^3}+\frac{x_1-m^i_1}{\sigma_1 r_i^2}\nonumber \\
&+2\frac{(x_2-m^i_2)^3(z^i_2-\hat\theta^i)}{\sigma_2(x_1-m^i_1)r_i^4}-2\frac{(x_2-m^i_2)(z^i_2-\hat\theta^i)}{\sigma_2(x_2-m^i_2)r_i^2)}+\frac{(x_2-m^i_2)^2}{\sigma_2r_i^4},\nonumber\\
\hat{H}^i_{12}=\hat{H}^i_{21}=&\frac{(x_1-m^i_1)(x_2-m^i_2)(r_i-z_1^i)}{\sigma_1r_i^3}+\frac{(x_1-m^i_1)(x_2-m^i_2)}{\sigma_1r_i^2}\nonumber \\
&+2\frac{(x_2-m^i_2)^2(z^i_2-\hat\theta^i)}{\sigma_2r_i^4}+\frac{(z^i_2-\hat\theta^i)}{\sigma_2r_i^2}-\frac{(x_1-m^i_1)(x_2-m^i_2)}{\sigma_2 r_i^4},\nonumber\\
\hat{H}^i_{13}=\hat{H}^i_{31}=&\frac{x_2-m^i_2}{\sigma_2 r_i^2},\nonumber\\
\hat{H}^i_{22}=& \frac{(x_2-m^i_2)^2(r_i-z^i_1)}{\sigma_1 r_i^3}+\frac{(r_i-z_1^i)}{\sigma_1 r_i}+\frac{(x_2-m^i_2)^2}{\sigma_1 r_i^2}\nonumber\\
&+2\frac{(x_1-m^i_1)(x_2-m^i_2)(z^i_2-\hat\theta^i)}{\sigma_2r^4}+\frac{(x_1-m^i_1)^2}{\sigma_2r^4},\nonumber\\
\hat{H}^i_{23}=\hat{H}^i_{32}=&\frac{x_1-m^i_1}{\sigma_2r_i^2}, \nonumber \\
\hat{H}^i_{33}=&\frac{1}{\sigma_2},\nonumber
\end{align}
and where
\begin{align}
r_i=&||x^n-m^i||,\nonumber \\
\hat\theta^i=&z^i_2-\left(\mbox{atan2}\left(\frac{m^i-x_2}{m^i_1-x_1}\right)-x_3+\frac{\pi}{2}\right).\nonumber
\end{align}

\bibliographystyle{plain}
\bibliography{Refs}
\end{document}